\documentclass{article}
\usepackage{amssymb}

\input{tcilatex}
\begin{document}

\begin{center}
\textbf{Existence local and global solution of multipoint Cauchy problem for
nonlocal nonlinear equations}

\textbf{Veli\ B. Shakhmurov}

Department of Mechanical Engineering, Okan University, Akfirat, Tuzla 34959
Istanbul, E-mail: veli.sahmurov@okan.edu.tr;

\textbf{Rishad Shahmurov}

shahmurov@hotmail.com

University of Alabama Tuscaloosa USA, AL 35487

A\textbf{bstract}
\end{center}

In this paper, the multipoint Cauchy problem for nonlocal nonlinear wave
type equat\i ons are studied.The equation involves a convolution integral
operator with a general kernel function whose Fourier transform is
nonnegative. We establish local and global existence and uniqueness of
solutions assuming enough smoothness on the initial data together with some
growth conditions on the nonlinear term


\begin{center}
\bigskip\ \ \textbf{AMS: 35Lxx, 35Qxx, 47D}

\textbf{1}. \textbf{Introduction}
\end{center}

The aim in this paper is to study the existence and uniqueness of solution
of the multipoint initial value problem (IVP) for nonlocal nonlinear wave
equat\i on 
\begin{equation}
u_{tt}-a\Delta u+b\ast u=\Delta \left[ g\ast f\left( u\right) \right] ,\text{
}x\in R^{n},\text{ }t\in \left( 0,\infty \right) ,  \tag{1.1}
\end{equation}%
\begin{equation}
u\left( x,0\right) =\varphi \left( x\right) +\dsum\limits_{k=1}^{m}\alpha
_{k}u\left( x,\lambda _{k}\right) ,\text{ for a.e. }x\in R^{n},  \tag{1.2}
\end{equation}%
\[
u_{t}\left( x,0\right) =\psi \left( x\right) +\dsum\limits_{k=1}^{m}\beta
_{k}u_{t}\left( x,\lambda _{k}\right) ,\text{ for a.e. }x\in R^{n}, 
\]%
where $m$ is an integer, $\lambda _{k}\in \left( 0,\infty \right) ,$ $\alpha
_{k},$ $\beta _{k}$ are complex numbers $g\left( x\right) $, $b\left(
x\right) $ are measurable functions on $\left( 0,\infty \right) $; $\ a\geq
0,$ $\Delta $ denotes the Laplace operator in $R^{n},$ $f(u)$ is the given
nonlinear function, $\varphi \left( x\right) $ and $\psi \left( x\right) $
are the given initial value functions. Note that for $\alpha _{k}=$ $\beta
_{k}=0$ we obtain classical Cauchy problem for nonlocal equation 
\begin{equation}
u_{tt}-a\Delta u+b\ast u=\Delta \left[ g\ast f\left( u\right) \right] ,\text{
}x\in R^{n},\text{ }t\in \left( 0,\infty \right) ,  \tag{1.3}
\end{equation}%
\[
u\left( 0,x\right) =\varphi \left( x\right) ,\text{ }u_{t}\left( 0,x\right)
=\psi \left( x\right) \text{ for a.e. }x\in R^{n}, 
\]

\bigskip The predictions of classical (local) elasticity theory become
inaccurate when the characteristic length of an elasticity problem is
comparable to the atomic length scale. To solution this situation, a
nonlocal theory of elasticity was introduced (see $[1-3]$ and the references
cited therein) and the main feature of the new theory is the fact that its
predictions were more down to earth than those of the classical theory. For
other generalizations of elasticity we refer the reader to $[4-6]$. \ The
global existence of the classical Cauchy problem for Boussinesq type
nonlocal equations has been studied by many authors (see $\left[ 7-11\right] 
$ ). Note that, the existence of solutions and regularity properties for
different type Boussinesq equations are considered e.g. in $\left[ \text{8-15%
}\right] $. Boussinesq type equations occur in a wide variety of physical
systems, such as in the propagation of longitudinal deformation waves in an
elastic rod, hydro-dynamical process in plasma, in materials science which
describe spinodal decomposition and in the absence of mechanical stresses
(see $\left[ 16-19\right] $).

The $L^{p}-$well-posedness of the class\i cal Cauchy problem $(1.3)$ depends
crucially on the presence of a suitable kernel. Then the question that
naturally arises is which of the possible forms of the kernel functions are
relevant for the global well-posedness of the multipoint initial-value
problem (IVP) $\left( 1.1\right) -\left( 1.2\right) $. In this study, as a
partial answer to this question, we consider multipoint IVP $(1.1)-\left(
1.2\right) $ with a general class of kernel functions and provide local,
global existence and blow-up results for the solutions of the problem $%
(1.1)-(1.2)$ in fram of $L^{p}$ spaces. The kernel functions most frequently
used in the literature are particular cases of this general class of kernel
functions. Note that nonlocal Cauchy problem for wave equations were studied
e.g. in $\left[ \text{20, 21}\right] .$

The strategy is to express the equation \ $\left( 1.1\right) $ as an
integral equation. To treat the nonlinearity as a small perturbation of the
linear part of the equation, the contraction mapping theorem is used. Also,
a priori estimates on $L^{p}$ norms of solutions of the linearized version
are utilized. The key step is the derivation of the uniform estimate of the
solutions of the linearized Boussinesq equation. The methods of harmonic
analysis, operator theory, interpolation of Banach Spaces and embedding
theorems in Sobolev spaces are the main tools implemented to carry out the
analysis.

In order to state our results precisely, we introduce some notations and
some function spaces.

\begin{center}
\textbf{Definitions and} \textbf{Background}
\end{center}

Let $E$ be a Banach space. $L^{p}\left( \Omega ;E\right) $ denotes the space
of strongly measurable $E$-valued functions that are defined on the
measurable subset $\Omega \subset R^{n}$ with the norm

\[
\left\Vert f\right\Vert _{p}=\left\Vert f\right\Vert _{L^{p}\left( \Omega
;E\right) }=\left( \int\limits_{\Omega }\left\Vert f\left( x\right)
\right\Vert _{E}^{p}dx\right) ^{\frac{1}{p}},\text{ }1\leq p<\infty ,\text{ }
\]

\[
\left\Vert f\right\Vert _{L^{\infty }\left( \Omega \right) }\ =\text{ess}%
\sup\limits_{x\in \Omega }\left\Vert f\left( x\right) \right\Vert _{E}. 
\]

Let $\mathbb{C}$ denote the set of complex numbers.\ For $E=\mathbb{C}$ the $%
L^{p}\left( \Omega ;E\right) $ denotes by $L^{p}\left( \Omega \right) .$

Let $E_{1}$ and $E_{2}$ be two Banach spaces. $\left( E_{1},E_{2}\right)
_{\theta ,p}$ for $\theta \in \left( 0,1\right) ,$ $p\in \left[ 1,\infty %
\right] $ denotes the interpolation spaces defined by $K$-method $\left[ 
\text{22, \S 1.3.2}\right] $.

\ Let $m$ be a positive integer. $W^{m,p}\left( \Omega \right) $ denotes the
Sobolev space, i.e. space of all functions $u\in L^{p}\left( \Omega \right) $
that have the generalized derivatives $\frac{\partial ^{m}u}{\partial
x_{k}^{m}}\in L^{p}\left( \Omega \right) ,$ $1\leq p\leq \infty $ with the
norm 
\[
\ \left\Vert u\right\Vert _{W^{m,p}\left( \Omega \right) }=\left\Vert
u\right\Vert _{L^{p}\left( \Omega \right) }+\sum\limits_{k=1}^{n}\left\Vert 
\frac{\partial ^{m}u}{\partial x_{k}^{m}}\right\Vert _{L^{p}\left( \Omega
\right) }<\infty . 
\]%
\ \ $\ \ $

Let $H^{s,p}\left( R^{n}\right) $, $\infty <s<\infty $ denotes fractionanal
Sobolev space of order $s$ which is defined as: 
\[
H^{s,p}=H^{s,p}\left( R^{n}\right) =\left( I-\Delta \right) ^{-\frac{s}{2}%
}L^{p}\left( R^{n}\right) 
\]%
with the norm 
\[
\left\Vert u\right\Vert _{H^{s,p}}=\left\Vert \left( I-\Delta \right) ^{%
\frac{s}{2}}u\right\Vert _{L^{p}\left( R^{n}\right) }<\infty . 
\]%
It clear that $H^{0,p}\left( R^{n}\right) =L^{p}\left( R^{n}\right) .$ It is
known that $H^{m,p}\left( R^{n}\right) =W^{m,p}\left( R^{n}\right) $ for the
positive integer $m$ (see e.g. $\left[ \text{23, \S\ 15}\right] $)$.$ For $%
p=2,$ the space $H^{s,p}\left( R^{n}\right) $ will be denoted by $%
H^{s}\left( R^{n}\right) .$ \ Let $S\left( R^{n}\right) $ denote Schwartz
class, i.e., the space of rapidly decreasing smooth functions on $R^{n},$
equipped with its usual topology generated by seminorms. Let $S^{^{\prime
}}\left( R^{n}\right) $ denote the space of all\ continuous linear operators 
$L:S\left( R^{n}\right) \rightarrow \mathbb{C},$ equipped with the bounded
convergence topology. Recall $S\left( R^{n}\right) $ is norm dense in $%
L^{p}\left( R^{n}\right) $ when $1\leq p<\infty .$

Let $1\leq p\leq q<\infty .$ \ A function $\Psi \in L_{\infty }(R^{n})$ is
called a Fourier multiplier from $L_{p}(R^{n})$ to $L_{q}(R^{n})$ if the map 
$B:$ $u\rightarrow F^{-1}\Psi (\xi )Fu$ for $u\in S(R^{n})$ is well defined
and extends to a bounded linear operator

\[
B:L_{p}(R^{n})\rightarrow L_{q}(R^{n}). 
\]
Let $L_{q}^{\ast }\left( E\right) $ denote the space of all $E-$valued
function space such that 
\[
\left\Vert u\right\Vert _{L_{q}^{\ast }\left( E\right) }=\left(
\int\limits_{0}^{\infty }\left\Vert u\left( t\right) \right\Vert _{E}^{q}%
\frac{dt}{t}\right) ^{\frac{1}{q}}<\infty ,\text{ }1\leq q<\infty ,\text{ }%
\left\Vert u\right\Vert _{L_{\infty }^{\ast }\left( E\right) }=\sup_{t\in
\left( 0,\infty \right) }\left\Vert u\left( t\right) \right\Vert _{E}. 
\]

Here, $F$ denote the Fourier transform. Fourier-analytic representation of
Besov spaces on $R^{n}$ is defined as:%
\[
B_{p,q}^{s}\left( R^{n}\right) =\left\{ u\in S^{^{\prime }}\left(
R^{n}\right) ,\right. \text{ } 
\]%
\[
\left\Vert u\right\Vert _{B_{p,q}^{s}\left( R^{n}\right) }=\left\Vert
F^{-1}t^{\varkappa -s}\left( 1+\left\vert \xi \right\vert ^{\frac{\varkappa 
}{2}}\right) e^{-t\left\vert \xi \right\vert ^{2}}Fu\right\Vert
_{L_{q}^{\ast }\left( L^{p}\left( R^{n}\right) \right) }\text{,} 
\]%
\[
\left\vert \xi \right\vert ^{2}=\dsum\limits_{k=1}^{n}\xi _{k}^{2}\text{, }%
\xi =\left( \xi _{1},\xi _{2},...,\xi _{n}\right) ,\left. p\in \left(
1,\infty \right) \text{, }q\in \left[ 1,\infty \right] \text{, }\varkappa
>s\right\} . 
\]%
\ It should be note that, the norm of Besov space does not depends on $%
\varkappa $ (see e.g. ( $\left[ \text{22, \S\ 2.3}\right] $). For $p=q$ the
space $B_{p,q}^{s}\left( R^{n}\right) $ will be denoted by $B_{p}^{s}\left(
R^{n}\right) .$

Sometimes we use one and the same symbol $C$ without distinction in order to
denote positive constants which may differ from each other even in a single
context. When we want to specify the dependence of such a constant on a
parameter, say $\alpha $, we write $C_{\alpha }$. Moreover, for $u$, $%
\upsilon >0$ the relation $u\lesssim \upsilon $ means that there exists a
constant $C>0$ \ independent on $u$ and $\upsilon $ such that 
\[
u\leq C\upsilon . 
\]

The paper is organized as follows: In Section 1, some definitions and
background are given. In Section 2, we obtain the existence of unique
solution and a priory estimates for solution of the linearized problem $%
(1.1)-\left( 1.2\right) .$ In Section 3, we show the existence and
uniqueness of local strong solution of the problem $(1.1)-\left( 1.2\right) $%
. In the Section 4 we show the same applications of the problem $%
(1.1)-\left( 1.2\right) .$

Sometimes we use one and the same symbol $C$ without distinction in order to
denote positive constants which may differ from each other even in a single
context. When we want to specify the dependence of such a constant on a
parameter, say $h$, we write $C_{h}$.

\begin{center}
\textbf{2. Estimates for linearized equation}
\end{center}

In this section, we make the necessary estimates for solutions of the Cauchy
problem%
\begin{equation}
u_{tt}-a\Delta u+b\ast u=g\left( x,t\right) ,\text{ }x\in R^{n},\text{ }t\in
\left( 0,T\right) ,\text{ }T\in \left( 0,\left. \infty \right] ,\right. 
\tag{2.1}
\end{equation}%
\begin{equation}
u\left( x,0\right) =\varphi \left( x\right) +\dsum\limits_{k=1}^{m}\alpha
_{k}u\left( x,\lambda _{k}\right) ,\text{ for a.e. }x\in R^{n},  \tag{2.2}
\end{equation}

\[
u_{t}\left( x,0\right) =\psi \left( x\right) +\dsum\limits_{k=1}^{m}\beta
_{k}u_{t}\left( x,\lambda _{k}\right) ,\text{ for a.e. }x\in R^{n}, 
\]%
Here, 
\[
X_{p}=L^{p}\left( R^{n}\right) \text{, }1\leq p\leq \infty ,\text{ }%
Y^{s,p}=H^{s,p}\left( R^{n}\right) ,\text{ }Y_{1}^{s,p}= 
\]

\[
H^{s,p}\left( R^{n}\right) \cap L^{1}\left( R^{n}\right) \text{, }Y_{\infty
}^{s,p}=H^{s,p}\left( R^{n}\right) \cap L^{\infty }\left( R^{n}\right) ,%
\text{ }Y_{p}^{s,p}=H^{s,p}\left( R^{n}\right) \cap L^{p}\left( R^{n}\right)
, 
\]%
\[
\left\Vert u\right\Vert _{Y_{p}^{s,2}}=\left\Vert u\right\Vert _{H^{s}\left(
R^{n}\right) }+\left\Vert u\right\Vert _{L^{p}\left( R^{n}\right) }<\infty
,1\leq p\leq \infty . 
\]

\bigskip Let $\hat{b}\left( \xi \right) $\ is the Fourier transformation of $%
b\left( x\right) ,$ i.e. $\hat{b}\left( \xi \right) =Fb$ and let%
\[
\eta =\eta \left( \xi \right) =\left[ a\left\vert \xi \right\vert ^{2}+\hat{b%
}\left( \xi \right) \right] ^{\frac{1}{2}}. 
\]
Here, 
\[
\varkappa _{k}=\varkappa _{k}\left( \xi \right) =\lambda _{k}\eta \left( \xi
\right) ,k=1,2,...,m. 
\]

\textbf{Condition 2.1. }Assume $a\geq 0,$ $b$ is an integrable function
whose $\hat{b}\left( \xi \right) \geq 0$ and $a+\hat{b}\left( \xi \right) >0 
$ for all $\ \xi \in R^{n}.$ Moreover, let \textbf{\ }%
\[
D_{0}\left( \xi \right) =1-\dsum\limits_{k=1}^{m}\left( \alpha _{k}+\beta
_{k}\right) \cos \varkappa _{k}+\dsum\limits_{i,j=1}^{m}\alpha _{i}\beta
_{j}\cos \left( \varkappa _{i}-\varkappa _{j}\right) \neq 0 
\]%
for all $\xi \in R^{n}$ with $\xi \neq 0.$

First we need the following lemmas:

\textbf{Lemma 2.1. }Let the Condition 2.1. holds. Then, the problem $\left(
2.1\right) -\left( 2.2\right) $ has a unique solution.

\textbf{Proof. }By using of the Fourier transform, we get from $(2.1)-\left(
2.2\right) $:%
\begin{equation}
\hat{u}_{tt}\left( \xi ,t\right) +\eta ^{2}\left( \xi \right) \hat{u}\left(
\xi ,t\right) =\hat{g}\left( \xi ,t\right) ,\text{ }  \tag{2.3}
\end{equation}%
\begin{equation}
\hat{u}\left( \xi ,0\right) =\hat{\varphi}\left( \xi \right)
+\dsum\limits_{k=1}^{m}\alpha _{k}\hat{u}\left( \xi ,\lambda _{k},\right) ,%
\text{ }\hat{u}_{t}\left( \xi ,0\right) =\hat{\psi}\left( \xi \right)
+\dsum\limits_{k=1}^{m}\beta _{k}\hat{u}\left( \xi ,\lambda _{k}\right) , 
\tag{2.4}
\end{equation}%
where $\hat{u}\left( \xi ,t\right) $ is a Fourier transform of $u\left(
x,t\right) $ with respect to $x$ and $\hat{\varphi}\left( \xi \right) ,$ $%
\hat{\psi}\left( \xi \right) $ are Fourier transform of $\varphi ,$ $\psi ,$
respectively.

By using the variation of constants it is easy to see that the general
solution of $\left( 2.3\right) $ is represented as%
\begin{equation}
\hat{u}\left( \xi ,t\right) =g_{1}\left( \xi \right) \cos \eta t+g_{2}\left(
\xi \right) \cos \eta t+\frac{1}{2\eta }\dint\limits_{0}^{t}\sin \eta \left(
t-\tau \right) \hat{g}\left( \xi ,\tau \right) d\tau ,  \tag{2.5}
\end{equation}%
where $g_{1}$, $g_{2}$ are general continuous differentiable functions. By
taking the multipoint condition $\left( 2.4\right) $ from $\left( 2.5\right) 
$ we get that $(2.3)-\left( 2.4\right) $ has a solution for $\xi \in R^{n},$
when $g_{1}$, $g_{2}$ are solution of the following system

\begin{equation}
g_{1}-\dsum\limits_{k=1}^{m}\alpha _{k}\gamma _{k}+g_{1}\cos \varkappa
_{k}+g_{2}\sin \varkappa _{k}=\hat{\varphi}\left( \xi \right) ,  \tag{2.6}
\end{equation}

\[
\eta g_{2}-\dsum\limits_{k=1}^{m}\beta _{k}\left[ \eta _{k}+\eta \left(
-g_{1}\sin \varkappa _{k}+g_{2}\cos \varkappa _{k}\right) \right] =\hat{\psi}%
\left( \xi \right) , 
\]%
where 
\[
\gamma _{k}=\gamma _{k}\left( \xi \right) =\frac{1}{2\eta }%
\dint\limits_{0}^{\lambda _{k}}\sin \eta \left( \lambda _{k}-\tau \right) 
\hat{g}\left( \xi ,\tau \right) d\tau ,\text{ } 
\]

\[
\mu _{k}=\mu _{k}\left( \xi \right) =-\frac{1}{2}\dint\limits_{0}^{\lambda
_{k}}\cos \eta \left( \lambda _{k}-\tau \right) \hat{g}\left( \xi ,\tau
\right) d\tau .\text{ } 
\]%
By Condition 2.1 we get 
\[
D\left( \xi \right) =\eta \left( \xi \right) D_{0}\left( \xi \right) \neq 0 
\]%
$\ $for all $\xi \neq 0.$ Solving the system $\left( 2.6\right) ,$ we obtain%
\begin{equation}
g_{1}=\frac{D_{1}\left( \xi \right) }{D\left( \xi \right) },\text{ }g_{2}=%
\frac{D_{1}\left( \xi \right) }{D\left( \xi \right) },  \tag{2.7}
\end{equation}%
where 
\[
D\left( \xi \right) =\left\vert 
\begin{array}{cc}
1-\dsum\limits_{k=1}^{m}\alpha _{k}\cos \varkappa _{k} & -\dsum%
\limits_{k=1}^{m}\alpha _{k}\sin \varkappa _{k} \\ 
--\eta \dsum\limits_{k=1}^{m}\beta _{k}\sin \varkappa _{k} & \eta \left(
1-\dsum\limits_{k=1}^{m}\beta _{k}\cos \varkappa _{k}\right)%
\end{array}%
\right\vert , 
\]%
\[
D_{1}\left( \xi \right) =\left\vert 
\begin{array}{cc}
\hat{\varphi}\left( \xi \right) +\dsum\limits_{k=1}^{m}\alpha _{k}\gamma _{k}
& -\dsum\limits_{k=1}^{m}\alpha _{k}\sin \varkappa _{k} \\ 
\hat{\psi}\left( \xi \right) +\dsum\limits_{k=1}^{m}\beta _{k}\mu _{k} & 
\eta \left( 1-\dsum\limits_{k=1}^{m}\beta _{k}\cos \varkappa _{k}\right)%
\end{array}%
\right\vert , 
\]%
\[
D_{2}\left( \xi \right) =\left\vert 
\begin{array}{cc}
1-\dsum\limits_{k=1}^{m}\alpha _{k}\cos \varkappa _{k} & \hat{\varphi}\left(
\xi \right) +\dsum\limits_{k=1}^{m}\alpha _{k}\gamma _{k} \\ 
-\eta \dsum\limits_{k=1}^{m}\beta _{k}\sin \varkappa _{k} & \hat{\psi}\left(
\xi \right) +\dsum\limits_{k=1}^{m}\beta _{k}\mu _{k}%
\end{array}%
\right\vert ; 
\]%
here, 
\[
\gamma _{k}\left( \xi \right) =\frac{1}{2\eta }\dint\limits_{0}^{\lambda
_{k}}\sin \eta \left( \lambda _{k}-\tau \right) \hat{g}\left( \xi ,\tau
\right) d\tau , 
\]%
\[
\mu _{k}\left( \xi \right) =-\frac{1}{2}\dint\limits_{0}^{\lambda _{k}}\cos
\eta \left( \lambda _{k}-\tau \right) \hat{g}\left( \xi ,\tau \right) d\tau
. 
\]%
Hense, problem $\left( 2.3\right) -\left( 2.4\right) $ has a unique solution
expressed as $\left( 2.5\right) ,$ where $g_{1}$ and $g_{2}$ are defined by $%
\left( 2.7\right) ,$ i.e. problem $(2.1)-\left( 2.2\right) $ has a unique
solution 
\begin{equation}
u\left( x,t\right) =F^{-1}\left[ C\left( \xi ,t\right) g_{1}\left( \xi
\right) \right] +\left[ F^{-1}S\left( \xi ,t\right) g_{2}\left( \xi \right) %
\right] +\frac{1}{2\eta }\dint\limits_{0}^{t}F^{-1}\left[ \sin \eta \left(
t-\tau \right) \hat{g}\left( \xi ,\tau \right) \right] d\tau .  \tag{2.8}
\end{equation}

\textbf{Theorem 2.1. }Let $\ $the Condition 2.1 holds and $s>\frac{n}{2}$.
Then for $\varphi ,$ $\psi ,$ $g\left( x,t\right) \in Y_{1}^{s,p}$ the
solution $\left( 2.1\right) -\left( 2.2\right) $ satisfies the following
uniformly in $t\in \left[ 0,T\right] $ estimate 
\begin{equation}
\left\Vert u\right\Vert _{X_{\infty }}+\left\Vert u_{t}\right\Vert
_{X_{\infty }}\leq C_{0}\left[ \left\Vert \varphi \right\Vert
_{Y_{1}^{s,p}}\right. +  \tag{2.9}
\end{equation}

\[
\left\Vert \psi \right\Vert _{Y_{1}^{s,p}}+\left. \dint\limits_{0}^{t}\left(
\left\Vert g\left( .,\tau \right) \right\Vert _{Y^{s,p}}+\left\Vert g\left(
.,\tau \right) \right\Vert _{X_{1}}\right) d\tau \right] , 
\]%
where the positive constant $C$ depends only on initial data.

\textbf{Proof. }From $\left( 2.7\right) $ we deduced that 
\[
g_{1}\left( \xi \right) =\eta ^{-1}\left( \xi \right) D_{0}^{-1}\left( \xi
\right) \left[ \eta \left( 1-\dsum\limits_{k=1}^{m}\beta _{k}\cos \varkappa
_{k}\right) \hat{\varphi}\left( \xi \right) +\right. 
\]%
\[
\hat{\psi}\left( \xi \right) \left( \dsum\limits_{k=1}^{m}\alpha _{k}\sin
\varkappa _{k}\right) +\left( \dsum\limits_{k=1}^{m}\beta _{k}\mu
_{k}\right) \left( \dsum\limits_{k=1}^{m}\alpha _{k}\sin \varkappa
_{k}\right) + 
\]%
\begin{equation}
\left. \eta \dsum\limits_{k=1}^{m}\alpha _{k}\gamma _{k}\left(
1-\dsum\limits_{k=1}^{m}\beta _{k}\cos \varkappa _{k}\right) \right] , 
\tag{2.10}
\end{equation}%
\[
\text{ }g_{2}\left( \xi \right) =\eta ^{-1}\left( \xi \right)
D_{0}^{-1}\left( \xi \right) \left[ \left( 1-\dsum\limits_{k=1}^{m}\alpha
_{k}\cos \varkappa _{k}\right) \hat{\psi}\left( \xi \right) +\right. 
\]%
\[
\dsum\limits_{k=1}^{m}\beta _{k}\mu _{k}\left(
1-\dsum\limits_{k=1}^{m}\alpha _{k}\cos \varkappa _{k}\right) +\hat{\varphi}%
\left( \xi \right) \dsum\limits_{k=1}^{m}\beta _{k}\sin \varkappa _{k}+ 
\]%
\[
\eta \left( \dsum\limits_{k=1}^{m}\beta _{k}\sin \varkappa _{k}\right)
\left( \dsum\limits_{k=1}^{m}\alpha _{k}\gamma _{k}\right) . 
\]%
Then, from $\left( 2.5\right) $, $\left( 2.8\right) $ and $\left(
2.10\right) $ we obtain that the solution $\left( 2.1\right) -\left(
2.2\right) $ can be expressed as 
\begin{equation}
u\left( x,t\right) =S_{1}\left( x,t\right) \varphi +S_{2}\left( x,t\right)
\psi +\Phi \left( g\right) \left( x,t\right) +\frac{1}{2\eta }%
\dint\limits_{0}^{t}F^{-1}\left[ \sin \eta \left( t-\tau \right) \hat{g}%
\left( \xi ,\tau \right) \right] d\tau ,  \tag{2.11}
\end{equation}%
where 
\[
S_{1}\left( x,t\right) \varphi =F^{-1}\left\{ D_{0}^{-1}\left( \xi \right) 
\left[ \left( 1-\dsum\limits_{k=1}^{m}\beta _{k}\cos \varkappa _{k}\right)
\sin \left( \eta t\right) \right] \right. ,+\text{ } 
\]%
\[
\left[ \eta ^{-1}\left( \xi \right) \dsum\limits_{k=1}^{m}\beta _{k}\sin
\varkappa _{k}\cos \left( \eta t\right) \right] \left. \hat{\varphi}\left(
\xi \right) \right\} , 
\]%
\[
S_{2}\left( x,t\right) \psi =F^{-1}\left\{ \left[ \eta ^{-1}\left( \xi
\right) D_{0}^{-1}\left( \xi \right) \left( \dsum\limits_{k=1}^{m}\alpha
_{k}\sin \varkappa _{k}\right) \sin \left( \eta t\right) \right] \right. +%
\text{ } 
\]

\[
\eta ^{-1}\left( \xi \right) D_{0}^{-1}\left( \xi \right) \left. \left[
\left( 1-\dsum\limits_{k=1}^{m}\alpha _{k}\cos \varkappa _{k}\right) \cos
\left( \eta t\right) \right] \right\} \hat{\psi}\left( \xi \right) , 
\]

\[
\Phi \left( x;t\right) =\Phi \left( g\right) \left( x;t\right) =\left[
\dsum\limits_{j=}^{4}\dsum\limits_{k=1}^{m}F^{-1}\Phi _{jk}\left( \xi
;t\right) \right] , 
\]%
where

\[
\Phi _{1k}\left( \xi ;t\right) =\frac{1}{2}D_{0}^{-1}\left( \xi \right)
\beta _{k}A_{1}\cos \left( \eta t\right) \dint\limits_{0}^{\lambda _{k}}\cos
\eta \left( \lambda _{k}-\tau \right) \hat{g}\left( \xi ,\tau \right) , 
\]%
\[
\Phi _{2k}\left( \xi ;t\right) =\frac{1}{2}D_{0}^{-1}\left( \xi \right)
\alpha _{k}A_{2}\cos \left( \eta t\right) \dint\limits_{0}^{\lambda
_{k}}\sin \eta \left( \lambda _{k}-\tau \right) \hat{g}\left( \xi ,\tau
\right) , 
\]%
\[
\Phi _{3k}\left( \xi ;t\right) =\frac{1}{2}D_{0}^{-1}\left( \xi \right)
\beta _{k}B_{1}\sin \left( \eta t\right) \dint\limits_{0}^{\lambda _{k}}\cos
\eta \left( \lambda _{k}-\tau \right) \hat{g}\left( \xi ,\tau \right) , 
\]%
\[
\Phi _{4k}\left( \xi ;t\right) =\frac{1}{2}D_{0}^{-1}\left( \xi \right)
\alpha _{k}B_{2}\sin \left( \eta t\right) \dint\limits_{0}^{\lambda
_{k}}\sin \eta \left( \lambda _{k}-\tau \right) \hat{g}\left( \xi ,\tau
\right) , 
\]%
here, 
\[
A_{1}=\dsum\limits_{k=1}^{m}\alpha _{k}\sin \varkappa _{k},\text{ }%
A_{2}=1-\dsum\limits_{k=1}^{m}\beta _{k}\cos \varkappa _{k},\text{ } 
\]

\[
B_{1}=\left( 1-\dsum\limits_{k=1}^{m}\alpha _{k}\cos \varkappa _{k}\right) ,%
\text{ \ }B_{2}=\dsum\limits_{k=1}^{m}\beta _{k}\cos \varkappa _{k}. 
\]%
By Condition 2.1, 
\[
D_{0}^{-1}\left( \xi \right) \eta ^{-1}\left( \xi \right) ,\text{ }\eta
^{-1}\left( \dsum\limits_{k=1}^{m}\alpha _{k}\sin \varkappa _{k}\right) ,%
\text{ }\eta ^{-1}\left( 1-\dsum\limits_{k=1}^{m}\alpha _{k}\cos \varkappa
_{k}\right) 
\]
are uniformly bounded. \ From $\left( 2.11\right) $ and $\left( 2.8\right) $
we obtain

\begin{equation}
\left\vert g_{1}\left( \xi \right) \right\vert \lesssim \left\vert \hat{%
\varphi}\left( \xi \right) \right\vert +\left\vert \hat{\psi}\left( \xi
\right) \right\vert +\left\vert \Phi \left( \xi \right) \right\vert , 
\tag{2.12}
\end{equation}%
\[
\left\vert g_{2}\left( \xi \right) \right\vert \lesssim \left\vert \hat{%
\varphi}\left( \xi \right) \right\vert +\left\vert \hat{\psi}\left( \xi
\right) \right\vert +\left\vert \Phi \left( \xi \right) \right\vert , 
\]%
where 
\[
\Phi \left( \xi \right) =\dsum\limits_{k=1}^{m}\dint\limits_{0}^{\lambda
_{k}}\hat{g}\left( \xi ,\tau \right) d\tau , 
\]%
Let $N\in \mathbb{N}$ and 
\[
\Pi _{N}=\left\{ \xi :\xi \in R^{n},\text{ }\left\vert \xi \right\vert \leq
N\right\} ,\text{ }\Pi _{N}^{\prime }=\left\{ \xi :\xi \in R^{n},\text{ }%
\left\vert \xi \right\vert \geq N\right\} . 
\]

From $\left( 2.8\right) $ we dedused that

\[
\left\Vert S_{1}\left( x,t\right) g_{1}\right\Vert _{X_{\infty }}+\left\Vert
S_{2}\left( x,t\right) g_{2}\right\Vert _{X_{\infty }}\lesssim 
\]%
\begin{equation}
\left\Vert F^{-1}C\left( \xi ,t\right) g_{1}\left( \xi \right) \right\Vert
_{L^{\infty }\left( \Pi _{N}\right) }+\left\Vert F^{-1}S\left( \xi ,t\right)
g_{2}\left( \xi \right) \right\Vert _{L^{\infty }\left( \Pi _{N}\right) }+ 
\tag{2.13}
\end{equation}%
\[
\left\Vert F^{-1}C\left( \xi ,t\right) g_{1}\left( \xi \right) \right\Vert
_{L^{\infty }\left( \Pi _{N}^{\prime }\right) }+\left\Vert F^{-1}S\left( \xi
,t\right) g_{2}\left( \xi \right) \right\Vert _{L^{\infty }\left( \Pi
_{N}^{\prime }\right) }. 
\]

\bigskip From $\left( 2.10\right) -$ $\left( 2.13\right) $ due to uniform
boundedness of $D_{0}^{-1}\left( \xi \right) $\ and $C\left( \xi ,t\right) $%
, $S\left( \xi ,t\right) $ we have 
\[
\left\Vert S_{1}\left( x,t\right) g_{1}\right\Vert _{X_{\infty }}+\left\Vert
S_{2}\left( x,t\right) g_{2}\right\Vert _{X_{\infty }}\lesssim \left\Vert
F^{-1}\hat{\varphi}\left( \xi \right) \right\Vert _{X_{\infty }}+\left\Vert
F^{-1}\hat{\psi}\left( \xi \right) \right\Vert _{X_{\infty }}+ 
\]%
\[
\left\Vert F^{-1}\dsum\limits_{k=1}^{m}\Phi \left( \xi \right) \right\Vert
_{X_{\infty }}, 
\]%
\[
\left\Vert S_{1}\left( x,t\right) g_{1}\right\Vert _{X_{\infty }}+\left\Vert
S_{2}\left( x,t\right) g_{2}\right\Vert _{X_{\infty }}\lesssim \left\Vert
F^{-1}\hat{\varphi}\left( \xi \right) \right\Vert _{X_{\infty }}+\left\Vert
F^{-1}\hat{\psi}\left( \xi \right) \right\Vert _{X_{\infty }}+ 
\]%
\[
\left\Vert F^{-1}\Phi \left( \xi \right) \right\Vert _{X_{\infty }}, 
\]%
In view of $\left( 2.12\right) ,$\ by using the Minkowski's inequality for
integrals from above\ we get 
\begin{equation}
\left\Vert F^{-1}C\left( \xi ,t\right) g_{1}\left( \xi \right) \right\Vert
_{L^{\infty }\left( \Pi _{N}\right) }+\left\Vert \left\Vert F^{-1}S\left(
\xi ,t\right) g_{2}\left( \xi \right) \right\Vert \right\Vert _{L^{\infty
}\left( \Pi _{N}\right) }\lesssim  \tag{2.14}
\end{equation}

\[
\left[ \left\Vert \varphi \right\Vert _{X_{1}}+\left\Vert \psi \right\Vert
_{X_{1}}+\left\Vert g\right\Vert _{X_{1}}\right] . 
\]%
Moreover, by $\left( 2.11\right) $ and $\left( 2.12\right) $ we have 
\[
\left\Vert F^{-1}C\left( \xi ,t\right) g_{1}\left( \xi \right) \right\Vert
_{L^{\infty }\left( \Pi _{N}^{\prime }\right) }+\left\Vert F^{-1}S\left( \xi
,t\right) g_{2}\left( \xi \right) \right\Vert _{L^{\infty }\left( \Pi
_{N}^{\prime }\right) }\lesssim 
\]%
\[
\left\Vert F^{-1}C\left( \xi ,t\right) \hat{\varphi}\left( \xi \right)
\right\Vert _{L^{\infty }\left( \Pi _{N}^{\prime }\right) }+\left\Vert
F^{-1}S\left( \xi ,t\right) \hat{\psi}\left( \xi \right) \right\Vert
_{L^{\infty }\left( \Pi _{N}^{\prime }\right) }+ 
\]%
\[
\left\Vert F^{-1}S\left( \xi ,t\right) \Phi \left( \xi \right) \right\Vert
_{L^{\infty }\left( \Pi _{N}^{\prime }\right) }\lesssim 
\]%
\begin{equation}
=\left\Vert F^{-1}\left( 1+\left\vert \xi \right\vert ^{2}\right) ^{-\frac{s%
}{2}}C\left( \xi ,t\right) \left( 1+\left\vert \xi \right\vert ^{2}\right) ^{%
\frac{s}{2}}\hat{\varphi}\left( \xi \right) \right\Vert _{L^{\infty }\left(
\Pi _{N}^{\prime }\right) }+  \tag{2.15}
\end{equation}%
\[
\left\Vert F^{-1}\left( 1+\left\vert \xi \right\vert ^{2}\right)
^{-s}S\left( \xi ,t\right) \left( 1+\left\vert \xi \right\vert \right) ^{%
\frac{s}{2}}\hat{\psi}\left( \xi \right) \right\Vert _{L^{\infty }\left( \Pi
_{N}^{\prime }\right) }+ 
\]%
\[
\left\Vert F^{-1}\left( 1+\left\vert \xi \right\vert ^{2}\right)
^{-s}S\left( \xi ,t\right) \left( 1+\left\vert \xi \right\vert \right) ^{%
\frac{s}{2}}\Phi \left( \xi \right) \right\Vert _{L^{\infty }\left( \Pi
_{N}^{\prime }\right) }+ 
\]%
By using $\left( 2.5\right) $, $\left( 2.7\right) $ and $\left( 2.12\right) $%
\ we get%
\[
\sup\limits_{\xi \in R^{n},t\in \left[ 0,T\right] }\left\vert \xi
\right\vert \left\vert ^{\left\vert \alpha \right\vert +\frac{n}{p}%
}D^{\alpha }\left[ \left( 1+\left\vert \xi \right\vert ^{2}\right) ^{-\frac{s%
}{2}}C\left( \xi ,t\right) \right] \right\vert \leq C_{2}, 
\]%
\ 
\begin{equation}
\sup\limits_{\xi \in R^{n},t\in \left[ 0,T\right] }\left\vert \xi
\right\vert \left\vert ^{\left\vert \alpha \right\vert +\frac{n}{p}%
}D^{\alpha }\left[ \left( 1+\left\vert \xi \right\vert ^{2}\right) ^{-\frac{s%
}{2}}S\left( \xi ,t\right) \right] \right\vert \leq C_{2}  \tag{2.16}
\end{equation}%
for $s>\frac{n}{p},$ $\alpha =\left( \alpha _{1},\alpha _{2},...,\alpha
_{n}\right) $, $\alpha _{k}\in \left\{ 0,1\right\} $, $\xi \in R^{n}$ and
uniformly in $t\in \left[ 0,T\right] .$ By multiplier theorems (see e.g. $%
\left[ 24\right] $) from $\left( 2.16\right) $ we get that the functions $%
\left( 1+\left\vert \xi \right\vert ^{2}\right) ^{-\frac{s}{2}}C\left( \xi
,t\right) ,$ $\left( 1+\left\vert \xi \right\vert ^{2}\right) ^{-\frac{s}{2}%
}S\left( \xi ,t\right) $ are $L^{p}\left( R^{n}\right) \rightarrow L^{\infty
}\left( R^{n}\right) $ Fourier multipliers. Then by Minkowski's inequality
for integrals, from $\left( 2.11\right) $ and $\left( 2.14\right) -\left(
2.16\right) $ we obtain%
\begin{equation}
\left\Vert F^{-1}C\left( \xi ,t\right) g_{1}\left( \xi \right) \right\Vert
_{L^{\infty }\left( \Pi _{N}^{\prime }\right) }+\left\Vert F^{-1}S\left( \xi
,t\right) g_{2}\left( \xi \right) \right\Vert _{L^{\infty }\left( \Pi
_{N}^{\prime }\right) }\lesssim  \tag{2.17}
\end{equation}

\[
\left[ \left\Vert \varphi \right\Vert _{Y^{s,p}}+\left\Vert \psi \right\Vert
_{Y^{s,p}}+\left\Vert g\right\Vert _{Y^{s,p}}\right] . 
\]

By reasoning as the above we have 
\begin{equation}
\left\Vert F^{-1}\Phi \left( \xi \right) \right\Vert _{X_{\infty }}\leq
C\dint\limits_{0}^{t}\left( \left\Vert g\left( .,\tau \right) \right\Vert
_{Y^{s}}+\left\Vert g\left( .,\tau \right) \right\Vert _{X_{1}}\right) d\tau
.  \tag{2.18}
\end{equation}

Thus, from $\left( 2.8\right) $ and $\left( 2.15\right) $ we obtain 
\begin{equation}
\left\Vert u\right\Vert _{X_{\infty }}\leq C\left[ \left\Vert \varphi
\right\Vert _{Y^{s,p}}+\left\Vert \varphi \right\Vert _{X_{1}}\right. + 
\tag{2.19}
\end{equation}

\[
\left\Vert \psi \right\Vert _{Y^{s,p}}+\left\Vert \psi \right\Vert
_{X_{1}}+\left. \dint\limits_{0}^{t}\left( \left\Vert g\left( .,\tau \right)
\right\Vert _{Y^{s,p}}+\left\Vert g\left( .,\tau \right) \right\Vert
_{X_{1}}\right) d\tau \right] . 
\]

By using $\left( 2.5\right) ,$ $\left( 2.7\right) $ and in view of $\left(
2.17\right) $ in similar way, we deduced the estimate of type $\left(
2.19\right) $ for $u_{t}$, i.e. we obtain the assertion.

\textbf{Theorem 2.2. }Let the Condition 2.1 holds and $s>\frac{n}{2}$. Then
for $\varphi ,$ $\psi ,$ $g\left( x,t\right) \in Y^{s,p}$ the solution of $%
\left( 2.1\right) -\left( 2.2\right) $ satisfies the following uniform
estimate%
\begin{equation}
\left( \left\Vert u\right\Vert _{Y^{s,p}}+\left\Vert u_{t}\right\Vert
_{Y^{s,p}}\right) \leq C_{0}\left( \left\Vert \varphi \right\Vert
_{Y^{s,p}}+\left\Vert \psi \right\Vert
_{Y^{s,p}}+\dint\limits_{0}^{t}\left\Vert g\left( .,\tau \right) \right\Vert
_{Y^{s,p}}d\tau \right) .  \tag{2.20}
\end{equation}

\textbf{Proof. }From $\left( 2.7\right) $ and $\left( 2.12\right) $ we have
the following uniform estimate 
\begin{equation}
\left( \left\Vert F^{-1}\left( 1+\left\vert \xi \right\vert ^{2}\right) ^{%
\frac{s}{2}}\hat{u}\right\Vert _{X_{p}}+\left\Vert F^{-1}\left( 1+\left\vert
\xi \right\vert ^{2}\right) ^{\frac{s}{2}}\hat{u}_{t}\right\Vert
_{X_{p}}\right) \leq  \tag{2.21}
\end{equation}

\[
C\left\{ \left\Vert F^{-1}\left( 1+\left\vert \xi \right\vert \right) ^{%
\frac{s}{2}}C\left( \xi ,t\right) \hat{\varphi}\right\Vert _{X_{p}}\right.
+\left\Vert F^{-1}\left( 1+\left\vert \xi \right\vert \right) ^{\frac{s}{2}%
}S\left( \xi ,t\right) \hat{\psi}\right\Vert _{X_{p}}+ 
\]

\[
\left. \dint\limits_{0}^{t}\left\Vert \left( 1+\left\vert \xi \right\vert
\right) ^{\frac{s}{2}}\hat{g}\left( .,\tau \right) \right\Vert _{X_{p}}d\tau
\right\} . 
\]

\bigskip By Condition 2.1 and by virtue of Fourier multiplier theorems (see $%
\left[ \text{24, \S\ 2.2\ }\right] $)\ we get that $C\left( \xi ,t\right) $, 
$S\left( \xi ,t\right) $ and $\Phi \left( \xi \right) $\ are Fourier
multipliers in $L^{p}\left( R^{n}\right) $ uniformly with respect to $t\in %
\left[ 0,T\right] .$ So, the estimate $\left( 2.21\right) $ by using the
Minkowski's inequality for integrals implies $\left( 2.20\right) .$

\begin{center}
\textbf{3. Local well posedness of IVP for nonlinear nonlocal equation}
\end{center}

In this section, we will show the local existence and uniqueness of solution
for the Cauchy problem $(1.1)-(1.2).$ For the study of the nonlinear problem 
$\left( 1.1\right) -\left( 1.2\right) $ we need the following lemmas

\textbf{Lemma 3.1} (Nirenberg's inequality) $\left[ 25\right] $. Assume that 
$u\in L^{p}\left( \Omega \right) $, $D^{m}u$ $\in L^{q}\left( \Omega \right) 
$, $p,q\in \left( 1,\infty \right) $. Then for $i$ with $0\leq i\leq m,$ $m>%
\frac{n}{q}$ we have 
\begin{equation}
\left\Vert D^{i}u\right\Vert _{r}\leq C\left\Vert u\right\Vert _{p}^{1-\mu
}\dsum\limits_{k=1}^{n}\left\Vert D_{k}^{m}u\right\Vert _{q}^{\mu }, 
\tag{3.1}
\end{equation}%
where%
\[
\frac{1}{r}=\frac{i}{m}+\mu \left( \frac{1}{q}-\frac{m}{n}\right) +\left(
1-\mu \right) \frac{1}{p},\text{ }\frac{i}{m}\leq \mu \leq 1. 
\]

\textbf{Lemma 3.2 }$\left[ 26\right] .$\textbf{\ }Assume that $u\in $ $%
W^{m,p}\left( \Omega \right) \cap L^{\infty }\left( \Omega \right) $ and $%
f\left( u\right) $ possesses continuous derivatives up to order $m\geq 1$.
Then $f\left( u\right) -f\left( 0\right) \in W^{m,p}\left( \Omega \right) $
and 
\[
\left\Vert f\left( u\right) -f\left( 0\right) \right\Vert _{p}\leq
\left\Vert f^{^{\left( 1\right) }}\left( u\right) \right\Vert _{\infty
}\left\Vert u\right\Vert _{p}, 
\]

\begin{equation}
\left\Vert D^{k}f\left( u\right) \right\Vert _{p}\leq
C_{0}\dsum\limits_{j=1}^{k}\left\Vert f^{\left( j\right) }\left( u\right)
\right\Vert _{\infty }\left\Vert u\right\Vert _{\infty }^{j-1}\left\Vert
D^{k}u\right\Vert _{p}\text{, }1\leq k\leq m,  \tag{3.2}
\end{equation}%
where $C_{0}$ $\geq 1$ is a constant.

Let%
\[
\text{ }X_{p}=L^{p}\left( R^{n}\right) ,\text{ }\left\Vert u\right\Vert
_{p}=\left\Vert u\right\Vert _{X_{p}},\text{ }Y=W^{2,p}\left( R^{n}\right) ,%
\text{ }E_{0}=\left( X_{p},Y\right) _{\frac{1}{2p},p}=B_{p}^{2-\frac{1}{p}%
}\left( R^{n}\right) . 
\]%
\textbf{Remark 3.1. }By using J.Lions-I. Petree result (see e.g $\left[ 
\text{21, \S\ 1.8.}\right] $) we obtain that the map $u\rightarrow u\left(
t_{0}\right) $, $t_{0}\in \left[ 0,T\right] $ is continuous and surjective
from $W^{2,p}\left( 0,T\right) $ onto $E_{0}$ and there is a constant $C_{1}$
such that 
\[
\left\Vert u\left( t_{0}\right) \right\Vert _{E_{0}}\leq C_{1}\left\Vert
u\right\Vert _{W^{2,p}\left( 0,T\right) },\text{ }1\leq p\leq \infty \text{.}
\]

First all of, we define the space $Y\left( T\right) =C\left( \left[ 0,T%
\right] ;Y_{\infty }^{2,p}\right) $ equipped with the norm defined by%
\[
\left\Vert u\right\Vert _{Y\left( T\right) }=\max\limits_{t\in \left[ 0,T%
\right] }\left\Vert u\right\Vert _{Y^{2,p}}+\max\limits_{t\in \left[ 0,T%
\right] }\left\Vert u\right\Vert _{X_{\infty }},\text{ }u\in Y\left(
T\right) . 
\]

It is easy to see that $Y\left( T\right) $ is a Banach space. For $\varphi $%
, $\psi \in Y^{2,p}$, let 
\[
M=\left\Vert \varphi \right\Vert _{Y^{2,p}}+\left\Vert \varphi \right\Vert
_{X_{\infty }}+\left\Vert \psi \right\Vert _{Y^{2,p}}+\left\Vert \psi
\right\Vert _{X_{\infty }}. 
\]

\textbf{Definition 3.1. }For any $T>0$ if $\upsilon ,$ $\psi \in Y_{\infty
}^{2,p}$ and $u$ $\in C\left( \left[ 0,T\right] ;Y_{\infty }^{2,p}\right) $
satisfies the equation $(1.1)-(1.2)$ then $u\left( x,t\right) $ is called
the continuous solution\ or the strong solution of the problem $(1.1)-(1.2).$
If $T<\infty $, then $u\left( x,t\right) $ is called the local strong
solution of the problem $(1.1)-(1.2).$ If $T=\infty $, then $u\left(
x,t\right) $ is called the global strong solution of the problem $%
(1.1)-(1.2) $.

\textbf{Condition 3.1. }Assume:

(1) Assume that the kernel $g$ is an integrable function whose Fourier
transform satisfies 
\[
0\leq \hat{g}\left( \xi \right) \lesssim \left( 1+\left\vert \xi \right\vert
^{2}\right) ^{-1}\text{ for all }\xi \in R^{n}; 
\]

(2) The Condition 2.1 holds, $\varphi ,$ $\psi $ $\in Y_{\infty }^{2,p}$ for 
$1<p<\infty $ $\ $and $\frac{n}{p}<2$;

(3) the function $u\rightarrow $ $f\left( x,t,u\right) $: $R^{n}\times \left[
0,T\right] \times E_{0}\rightarrow E$ is a measurable in $\left( x,t\right)
\in R^{n}\times \left[ 0,T\right] $ for $u\in E_{0};$ $f\left( x,t,u\right) $%
. Moreover, $F\left( x,t,u\right) $ is continuous in $u\in E_{0}$ and $%
f\left( x,t,u\right) \in C^{\left( 3\right) }\left( E_{0};E\right) $
uniformly with respect to $x\in R^{n},$ $t\in \left[ 0,T\right] .$

Main aim of this section is to prove the following result:

\textbf{Theorem 3.1. }Let the Condition 3.1. holds. Then problem $\left(
1.1\right) -\left( 1.2\right) $ has a unique local strange solution $u\in
C^{\left( 2\right) }\left( \left[ 0,\right. \left. T_{0}\right) ;Y_{\infty
}^{2,p}\right) $, where $T_{0}$ is a maximal time interval that is
appropriately small relative to $M$. Moreover, if

\begin{equation}
\sup_{t\in \left[ 0\right. ,\left. T_{0}\right) }\left( \left\Vert
u\right\Vert _{Y_{\infty }^{2,p}}+\left\Vert u_{t}\right\Vert _{Y_{\infty
}^{2,p}}\right) <\infty  \tag{3.3}
\end{equation}%
then $T_{0}=\infty .$

\textbf{Proof. }First, we are going to prove the existence and the
uniqueness of the local continuous solution of the problem $(1.1)-\left(
1.2\right) $ by contraction mapping principle. Consider a map $G$ on $%
Y\left( T\right) $ such that $G(u)$ is the solution of the Cauchy problem%
\begin{equation}
G_{tt}\left( u\right) -a\Delta G\left( u\right) =\Delta \left[ g\ast f\left(
G\left( u\right) \right) \right] ,\text{ }x\in R^{n},\text{ }t\in \left(
0,T\right) ,  \tag{3.4}
\end{equation}%
\[
G\left( u\right) \left( x,0\right) =\varphi \left( x\right)
+\dsum\limits_{k=1}^{m}\alpha _{k}G\left( u\right) \left( x,\lambda
_{k}\right) ,\text{ for a.e. }x\in R^{n}, 
\]%
\[
G_{t}\left( u\right) \left( 0,x\right) =\psi \left( x\right)
+\dsum\limits_{k=1}^{m}\beta _{k}G_{t}\left( u\right) \left( x,\lambda
_{k}\right) ,\text{ for a.e. }x\in R^{n}. 
\]

From Lemma 3.2 we know that $F(u)\in $ $L^{p}\left( 0,T;Y_{\infty
}^{2,p}\right) $ for any $T>0$. Thus, by Lemma 2.1, problem $\left(
3.4\right) $ has a solution which can be written as%
\begin{equation}
G\left( u\right) \left( x,t\right) =\left[ S_{1}\left( x,t\right) \varphi
+S_{2}\left( x,t\right) \psi +\Phi \left( g\ast f\left( G\left( u\right)
\right) \right) \right] +  \tag{3.5}
\end{equation}%
\[
\frac{1}{2\eta }\dint\limits_{0}^{t}F^{-1}\left[ \sin \eta \left( t-\tau
\right) \left\vert \xi \right\vert ^{2}\hat{g}\left( \xi \right) \hat{f}%
\left( G\left( u\right) \left( \xi \right) \right) \right] d\tau , 
\]%
where $S_{1}\left( x,t\right) $, $S_{2}\left( x,t\right) ,$ $\Phi $ are
operator functions defined by $\left( 2.10\right) $ and $\left( 2.11\right) $%
, $\ $where $g$ replaced by $g\ast f\left( G\left( u\right) \right) .$ From
Lemma 3.2 it is easy to see that the map $G$ is well defined for $f\in
C^{\left( 2\right) }\left( X_{0};\mathbb{C}\right) $. We put 
\[
Q\left( M;T\right) =\left\{ u\mid u\in Y\left( T\right) \text{, }\left\Vert
u\right\Vert _{Y\left( T\right) }\leq M+1\right\} . 
\]

First, by reasoning as in $\left[ 9\right] $\ let us prove that the map $G$
has a unique fixed point in $Q\left( M;T\right) .$ For this aim, it is
sufficient to show that the operator $G$ maps $Q\left( M;T\right) $ into $%
Q\left( M;T\right) $ and $G:$ $Q\left( M;T\right) $ $\rightarrow $ $Q\left(
M;T\right) $ is strictly contractive if $T$ is appropriately small relative
to $M.$ Consider the function \ $\bar{f}\left( \xi \right) $: $\left[
0,\right. $ $\left. \infty \right) \rightarrow \left[ 0,\right. $ $\left.
\infty \right) $ defined by 
\[
\ \bar{f}\left( \xi \right) =\max\limits_{\left\vert x\right\vert \leq \xi
}\left\{ \left\Vert f^{\left( 1\right) }\left( x\right) \right\Vert _{%
\mathbb{C}},\left\Vert f^{\left( 2\right) }\left( x\right) \right\Vert _{%
\mathbb{C}}\text{ }\right\} ,\text{ }\xi \geq 0. 
\]

It is clear to see that the function $\bar{f}\left( \xi \right) $ is
continuous and nondecreasing on $\left[ 0,\right. $ $\left. \infty \right) .$
From Lemma 3.2 we have\qquad

\[
\left\Vert f\left( u\right) \right\Vert _{Y^{2,p}}\leq \left\Vert f^{\left(
1\right) }\left( u\right) \right\Vert _{X_{\infty }}\left\Vert u\right\Vert
_{X_{p}}+\left\Vert f^{\left( 1\right) }\left( u\right) \right\Vert
_{X_{\infty }}\left\Vert Du\right\Vert _{X_{p}}+ 
\]

\begin{equation}
C_{0}\left[ \left\Vert f^{\left( 1\right) }\left( u\right) \right\Vert
_{X_{\infty }}\left\Vert u\right\Vert _{X_{p}}+\left\Vert f^{\left( 2\right)
}\left( u\right) \right\Vert _{X_{\infty }}\left\Vert u\right\Vert
_{X_{\infty }}\left\Vert D^{2}u\right\Vert _{X_{p}}\right] \leq  \tag{3.6}
\end{equation}

\ 
\[
2C_{0}\bar{f}\left( M+1\right) \left( M+1\right) \left\Vert u\right\Vert
_{Y^{2,p}}. 
\]%
In view of the assumpt\i on (1) and by using Minkowski's inequality for
integrals\"{o} we obtain from $\left( 3.5\right) $:%
\begin{equation}
\left\Vert G\left( u\right) \right\Vert _{X_{\infty }}\lesssim \left\Vert
\varphi \right\Vert _{\infty }+\left\Vert \psi \right\Vert _{\infty
}+\dint\limits_{0}^{t}\left\Vert \Delta \left[ g\ast f\left( G\left(
u\right) \right) \right] \left( x,\tau \right) d\tau \right\Vert _{\infty },
\tag{3.7}
\end{equation}%
\begin{equation}
\left\Vert G\left( u\right) \right\Vert _{Y^{2,p}}\lesssim \left\Vert
\varphi \right\Vert _{Y^{2,p}}+\left\Vert \psi \right\Vert
_{Y^{2,p}}+\dint\limits_{0}^{t}\left\Vert \Delta \left[ g\ast f\left(
G\left( u\right) \right) \right] \left( x,\tau \right) d\tau \right\Vert
_{Y^{2,p}}d\tau .  \tag{3.8}
\end{equation}%
Thus, from $\left( 3.6\right) -\left( 3.8\right) $ and Lemma 3.2 we get 
\[
\left\Vert G\left( u\right) \right\Vert _{Y\left( T\right) }\leq M+T\left(
M+1\right) \left[ 1+2C_{0}\left( M+1\right) \bar{f}\left( M+1\right) \right]
. 
\]%
If $T$ satisfies 
\begin{equation}
T\leq \left\{ \left( M+1\right) \left[ 1+2C_{0}\left( M+1\right) \bar{f}%
\left( M+1\right) \right] \right\} ^{-1},  \tag{3.9}
\end{equation}%
then 
\[
\left\Vert Gu\right\Vert _{Y\left( T\right) }\leq M+1. 
\]%
Therefore, if $\left( 3.9\right) $ holds, then $G$ maps $Q\left( M;T\right) $
into $Q\left( M;T\right) .$ Now, we are going to prove that the map $G$ is
strictly contractive. Assume $T>0$ and $u_{1},$ $u_{2}\in $ $Q\left(
M;T\right) $ given. We get%
\[
G\left( u_{1}\right) -G\left( u_{2}\right) = 
\]%
\[
\dint\limits_{0}^{t}F^{-1}S\left( t-\tau ,\xi \right) \left\vert \xi
\right\vert ^{2}\hat{g}\left( \xi \right) \left[ \hat{f}\left( u_{1}\right)
\left( \xi ,\tau \right) -\hat{f}\left( u_{2}\right) \left( \xi ,\tau
\right) \right] d\tau ,\text{ }t\in \left( 0,T\right) . 
\]%
By using the assumption (3) and the mean value theorem, we obtain%
\[
\hat{f}\left( u_{1}\right) -\hat{f}\left( u_{2}\right) =\hat{f}^{\left(
1\right) }\left( u_{2}+\eta _{1}\left( u_{1}-u_{2}\right) \right) \left(
u_{1}-u_{2}\right) ,\text{ } 
\]

\[
D_{\xi }\left[ \hat{f}\left( u_{1}\right) -\hat{f}\left( u_{2}\right) \right]
=\hat{f}^{\left( 2\right) }\left( u_{2}+\eta _{2}\left( u_{1}-u_{2}\right)
\right) \left( u_{1}-u_{2}\right) D_{\xi }u_{1}+\text{ } 
\]%
\[
\hat{f}^{\left( 1\right) }\left( u_{2}\right) \left( D_{\xi }u_{1}-D_{\xi
}u_{2}\right) , 
\]%
\[
D_{\xi }^{2}\left[ \hat{f}\left( u_{1}\right) -\hat{f}\left( u_{2}\right) %
\right] =\hat{f}^{\left( 3\right) }\left( u_{2}+\eta _{3}\left(
u_{1}-u_{2}\right) \right) \left( u_{1}-u_{2}\right) \left( D_{\xi
}u_{1}\right) ^{2}+\text{ } 
\]%
\[
\hat{f}^{\left( 2\right) }\left( u_{2}\right) \left( D_{\xi }u_{1}-D_{\xi
}u_{2}\right) \left( D_{\xi }u_{1}+D_{\xi }u_{2}\right) + 
\]%
\[
\hat{f}^{\left( 2\right) }\left( u_{2}+\eta _{4}\left( u_{1}-u_{2}\right)
\right) \left( u_{1}-u_{2}\right) D_{\xi }^{2}u_{1}+\hat{f}^{\left( 1\right)
}\left( u_{2}\right) \left( D_{\xi }^{2}u_{1}-D_{\xi }^{2}u_{2}\right) , 
\]%
where $0<\eta _{i}<1,$ $i=1,2,3,4.$ Thus, using H\"{o}lder's and Nirenberg's
inequality, we have%
\begin{equation}
\left\Vert \hat{f}\left( u_{1}\right) -\hat{f}\left( u_{2}\right)
\right\Vert _{X_{\infty }}\leq \bar{f}\left( M+1\right) \left\Vert
u_{1}-u_{2}\right\Vert _{X_{\infty }},  \tag{3.10}
\end{equation}%
\begin{equation}
\left\Vert \hat{f}\left( u_{1}\right) -\hat{f}\left( u_{2}\right)
\right\Vert _{X_{p}}\leq \bar{f}\left( M+1\right) \left\Vert
u_{1}-u_{2}\right\Vert _{X_{p}},  \tag{3.11}
\end{equation}%
\begin{equation}
\left\Vert D_{\xi }\left[ \hat{f}\left( u_{1}\right) -\hat{f}\left(
u_{2}\right) \right] \right\Vert _{X_{p}}\leq \left( M+1\right) \bar{f}%
\left( M+1\right) \left\Vert u_{1}-u_{2}\right\Vert _{X_{\infty }}+ 
\tag{3.12}
\end{equation}%
\[
\bar{f}\left( M+1\right) \left\Vert \hat{f}\left( u_{1}\right) -\hat{f}%
\left( u_{2}\right) \right\Vert _{X_{p}}, 
\]%
\[
\left\Vert D_{\xi }^{2}\left[ \hat{f}\left( u_{1}\right) -\hat{f}\left(
u_{2}\right) \right] \right\Vert _{X_{p}}\leq \left( M+1\right) \bar{f}%
\left( M+1\right) \left\Vert u_{1}-u_{2}\right\Vert _{X_{\infty }}\left\Vert
D_{\xi }^{2}u_{1}\right\Vert _{Y^{2,p}}^{2}+ 
\]%
\[
\bar{f}\left( M+1\right) \left\Vert D_{\xi }\left( u_{1}-u_{2}\right)
\right\Vert _{Y^{2,p}}\left\Vert D_{\xi }\left( u_{1}+u_{2}\right)
\right\Vert _{Y^{2,p}}+ 
\]%
\[
\bar{f}\left( M+1\right) \left\Vert u_{1}-u_{2}\right\Vert _{X_{\infty
}}\left\Vert D_{\xi }^{2}u_{1}\right\Vert _{X_{p}}+\bar{f}\left( M+1\right)
\left\Vert D_{\xi }\left( u_{1}-u_{2}\right) \right\Vert _{X_{p}}\leq 
\]%
\begin{equation}
C^{2}\bar{f}\left( M+1\right) \left\Vert u_{1}-u_{2}\right\Vert _{X_{\infty
}}\left\Vert u_{1}\right\Vert _{X_{\infty }}\left\Vert D_{\xi
}^{2}u_{1}\right\Vert _{X_{p}}+  \tag{3.13}
\end{equation}%
\[
C^{2}\bar{f}\left( M+1\right) \left\Vert u_{1}-u_{2}\right\Vert _{X_{\infty
}}^{\frac{1}{2}}\left\Vert D_{\xi }^{2}\left( u_{1}-u_{2}\right) \right\Vert
_{X_{p}}\left\Vert u_{1}+u_{2}\right\Vert _{X_{\infty }}^{\frac{1}{2}%
}\left\Vert D_{\xi }^{2}\left( u_{1}+u_{2}\right) \right\Vert _{X_{p}} 
\]%
\[
+\left( M+1\right) \bar{f}\left( M+1\right) \left\Vert
u_{1}-u_{2}\right\Vert _{X_{\infty }}+\bar{f}\left( M+1\right) \left\Vert
D_{\xi }^{2}\left( u_{1}-u_{2}\right) \right\Vert _{X_{p}}\leq 
\]%
\[
3C^{2}\left( M+1\right) ^{2}\bar{f}\left( M+1\right) \left\Vert
u_{1}-u_{2}\right\Vert _{X_{\infty }}+2C^{2}\left( M+1\right) \bar{f}\left(
M+1\right) \left\Vert D_{\xi }^{2}\left( u_{1}-u_{2}\right) \right\Vert
_{X_{p}}, 
\]%
where $C$ is the constant in Lemma $3.1$. From $\left( 3.10\right) -\left(
3.11\right) $, using Minkowski's inequality for integrals, Fourier
multiplier theorems in $X_{p}$ spaces and Young's inequality, we obtain%
\[
\left\Vert G\left( u_{1}\right) -G\left( u_{2}\right) \right\Vert _{Y\left(
T\right) }\leq \dint\limits_{0}^{t}\left\Vert u_{1}-u_{2}\right\Vert
_{X_{\infty }}d\tau +\dint\limits_{0}^{t}\left\Vert u_{1}-u_{2}\right\Vert
_{Y^{2,p}}d\tau + 
\]%
\[
\dint\limits_{0}^{t}\left\Vert f\left( u_{1}\right) -f\left( u_{2}\right)
\right\Vert _{X_{\infty }}d\tau +\dint\limits_{0}^{t}\left\Vert f\left(
u_{1}\right) -f\left( u_{2}\right) \right\Vert _{Y^{2,p}}d\tau \leq 
\]%
\[
T\left[ 1+C_{1}\left( M+1\right) ^{2}\bar{f}\left( M+1\right) \right]
\left\Vert u_{1}-u_{2}\right\Vert _{Y\left( T\right) }, 
\]%
where $C_{1}$ is a constant. If $T$ satisfies $\left( 3.9\right) $ and the
following inequality holds 
\begin{equation}
T\leq \frac{1}{2}\left[ 1+C_{1}\left( M+1\right) ^{2}\bar{f}\left(
M+1\right) \right] ^{-1},  \tag{3.14}
\end{equation}%
then 
\[
\left\Vert Gu_{1}-Gu_{2}\right\Vert _{Y\left( T\right) }\leq \frac{1}{2}%
\left\Vert u_{1}-u_{2}\right\Vert _{Y\left( T\right) }. 
\]

That is, $G$ is a contructive map. By contraction mapping principle we know
that $G(u)$ has a fixed point $u(x,t)\in $ $Q\left( M;T\right) $ that is a
solution of $(1.1)-(1.2)$. From $\left( 2.9\right) -\left( 2.11\right) $ we
get that $u$ is a solution of the following integral equation 
\[
u\left( x,t\right) =S_{1}\left( t\right) g_{1}+S_{2}\left( t\right) g_{2}- 
\]%
\[
\dint\limits_{0}^{t}F^{-1}\left[ S\left( t-\tau ,\xi \right) \left\vert \xi
\right\vert ^{2}\hat{g}\left( \xi \right) \hat{f}\left( u\right) \left( \xi
,\tau \right) \right] d\tau ,\text{ }t\in \left( 0,T\right) . 
\]

Let us show that this solution is a unique in $Y\left( T\right) $. Let $%
u_{1} $, $u_{2}\in Y\left( T\right) $ are two solution of the problem $%
(1.1)-(1.2)$. Then%
\begin{equation}
u_{1}-u_{2}=\dint\limits_{0}^{t}F^{-1}\left[ S\left( t-\tau ,\xi \right)
\left\vert \xi \right\vert ^{2}\hat{g}\left( \xi \right) \left( \hat{f}%
\left( u_{1}\right) \left( \xi ,\tau \right) -\hat{f}\left( u_{2}\right)
\left( \xi ,\tau \right) \right) \right] d\tau .  \tag{3.15}
\end{equation}%
By the definition of the space $Y\left( T\right) $, we can assume that%
\[
\left\Vert u_{1}\right\Vert _{X_{\infty }}\leq C_{1}\left( T\right) ,\text{ }%
\left\Vert u_{1}\right\Vert _{X_{\infty }}\leq C_{1}\left( T\right) . 
\]

Hence, by Minkowski's inequality for integrals and Theorem 2.2 we obtain
from $\left( 3.15\right) $

\begin{equation}
\left\Vert u_{1}-u_{2}\right\Vert _{Y^{2,p}}\leq C_{2}\left( T\right) \text{ 
}\dint\limits_{0}^{t}\left\Vert u_{1}-u_{2}\right\Vert _{Y^{2,p}}d\tau . 
\tag{3.16}
\end{equation}

From $(3.16)$ and Gronwall's inequality, we have $\left\Vert
u_{1}-u_{2}\right\Vert _{Y^{2,p}}=0$, i.e. problem $(1.1)-(1.2)$ has a
unique solution which belongs to $Y\left( T\right) .$ That is, we obtain the
first part of the assertion.

Now, let $\left[ 0\right. ,\left. T_{0}\right) $ be the maximal time
interval of existence for $u\in Y\left( T_{0}\right) $. It remains only to
show that if $(3.3)$ is satisfied, then $T_{0}=\infty $. Assume contrary
that, $\left( 3.3\right) $ holds and $T_{0}<\infty .$ For $T\in \left[
0\right. ,\left. T_{0}\right) ,$ we consider the following integral equation

\begin{equation}
\upsilon \left( x,t\right) =S_{1}\left( t\right) u\left( x,T\right)
+S_{2}\left( t\right) u_{t}\left( x,T\right) -  \tag{3.17}
\end{equation}

\[
\dint\limits_{0}^{t}F^{-1}\left[ S\left( t-\tau ,\xi \right) \left\vert \xi
\right\vert ^{2}\hat{g}\left( \xi \right) \hat{f}\left( \upsilon \right)
\left( \xi ,\tau \right) \right] d\tau ,\text{ }t\in \left( 0,T\right) . 
\]%
By virtue of $(3.3)$, for $T^{\prime }>T$ we have 
\[
\sup_{t\in \left[ 0\right. ,\left. T\right) }\left( \left\Vert u\right\Vert
_{Y^{2,p}}+\left\Vert u\right\Vert _{X_{\infty }}+\left\Vert
u_{t}\right\Vert _{Y^{2,p}}+\left\Vert u_{t}\right\Vert _{X_{\infty
}}\right) <\infty . 
\]

By reasoning as a first part of theorem and by contraction mapping
principle, there is a $T^{\ast }\in \left( 0,T_{0}\right) $ such that for
each $T\in \left[ 0\right. ,\left. T_{0}\right) ,$ the equation $\left(
3.17\right) $ has a unique solution $\upsilon \in Y\left( T^{\ast }\right) .$
The estimates $\left( 3.9\right) $ and $\left( 3.14\right) $ imply that $%
T^{\ast }$ can be selected independently of $T\in \left[ 0\right. ,\left.
T_{0}\right) .$ Set $T=T_{0}-\frac{T^{\ast }}{2}$ and define 
\begin{equation}
\tilde{u}\left( x,t\right) =\left\{ 
\begin{array}{c}
u\left( x,t\right) ,\text{ }t\in \left[ 0,T\right] \\ 
\upsilon \left( x,t-T\right) \text{, }t\in \left[ T,T_{0}+\frac{T^{\ast }}{2}%
\right]%
\end{array}%
\right. .  \tag{3.18}
\end{equation}

By construction $\tilde{u}\left( x,t\right) $ is a solution of the problem $%
(1.1)-(1.2)$ on $\left[ T,T_{0}+\frac{T^{\ast }}{2}\right] $ and in view of
local uniqueness, $\tilde{u}\left( x,t\right) $ extends $u.$ This is against
to the maximality of $\left[ 0\right. ,\left. T_{0}\right) $, i.e we obtain $%
T_{0}=\infty .$

Consider the problem $\left( 1.1\right) -\left( 1.2\right) ,$ when $\varphi
, $ $\psi \in H^{s}.$

We first need two lemmas concerning the behaviour of the nonlinear term $[$%
8, 13, 27$]$.

\textbf{\ Lemma 3.3.} Let $s\geq 0,$ $f\in C^{\left[ s\right] +1}\left(
R\right) $ with $f(0)=0$. Then for any $u\in H^{s}\cap L^{\infty }$, we have 
$f(u)\in H^{s}\cap L^{\infty }.$ Moreover there is some constant $A(M)$
depending on $M$ such that for all $u\in H^{s}\cap L^{\infty }$ with $%
\left\Vert u\right\Vert _{L^{m}}\leq M,$%
\[
\left\Vert f(u)\right\Vert _{H^{s}}\leq A\left( M\right) \left\Vert
u)\right\Vert _{H^{s}}. 
\]%
\textbf{Lemma 3.4. } Let $s\geq 0,$ $f\in C^{\left[ s\right] +1}\left(
R\right) $. Then for for any $M$ there is some constant $B(M)$ depending on $%
M$ such that for all $u$, $\upsilon \in H^{s}\cap L^{\infty }$ with $%
\left\Vert u\right\Vert _{L^{m}}\leq M,$ $\left\Vert \upsilon \right\Vert
_{L^{\infty }}\leq M,$ $\left\Vert u\right\Vert _{H^{8}}\leq M,$ $\left\Vert
\upsilon \right\Vert _{H^{s}}\leq M,$%
\[
\left\Vert f(u)-f(\upsilon \right\Vert _{H^{s}}\leq B\left( M\right)
\left\Vert u-\upsilon \right\Vert _{H^{s}},\text{ }\left\Vert
f(u)-f(\upsilon \right\Vert _{L^{\infty }}\leq B\left( M\right) \left\Vert
u-\upsilon \right\Vert _{L^{\infty }}. 
\]

By reasoning as in $\left[ \text{13, Lemma 3.4}\right] $ we have

\textbf{Corollary 3.1.} Let $s>\frac{n}{2},$ $f\in C^{\left[ s\right]
+1}\left( R\right) $. Then for for any $B$ there is some constant $B(M)$
depending on $M$ such that for all $u$, $\upsilon \in H^{s}$ with $%
\left\Vert u\right\Vert _{H^{s}}\leq M,$ $\left\Vert \upsilon \right\Vert
_{H^{s}}\leq M,$%
\[
\left\Vert f(u)-f(\upsilon \right\Vert _{H^{s}}\leq B\left( M\right)
\left\Vert u-\upsilon \right\Vert _{H^{s}}. 
\]

\bigskip \textbf{Lemma 3.5. } If $s>0$, then $Y_{\infty }^{s,2}$ is an
algebra. Moreover, for \ $f,$ $g\in Y_{\infty }^{s,2},$ 
\[
\left\Vert fg\right\Vert _{H^{s}}\leq C\left[ \left\Vert f\right\Vert
_{\infty }+\left\Vert g\right\Vert _{H^{s}}+\left\Vert f\right\Vert
_{H^{s}}+\left\Vert g\right\Vert _{\infty }\right] . 
\]

\textbf{\ Lemma 3.6 }$\left[ \text{28, Lemma X 4}\right] $\textbf{.} Let $%
s\geq 0,$ $f\in C^{\left[ s\right] +1}\left( R\right) $ and $f\left(
u\right) =O\left( \left\vert u\right\vert ^{\alpha +1}\right) $ \ for $%
u\rightarrow 0$, $\alpha \geq 1$ be a positive integer. If $u\in Y_{\infty
}^{s,2}$ and $\left\Vert u\right\Vert _{\infty }\leq M$, \ then 
\[
\left\Vert f(u)\right\Vert _{H^{s}}\leq C\left( M\right) \left[ \left\Vert
u\right\Vert _{H^{s}}\left\Vert u\right\Vert _{\infty }^{\alpha }\right] , 
\]%
\[
\left\Vert f(u)\right\Vert _{1}\leq C\left( M\right) \left\Vert u\right\Vert
_{2}^{2}\left\Vert u\right\Vert _{\infty }^{\alpha -1}. 
\]

\textbf{Lemma 3.7 }$\left[ \text{13, Lemma 3.4}\right] $\textbf{.} Let $%
s\geq 0,$ $f\in C^{\left[ s\right] +1}\left( R\right) $ and $f\left(
u\right) =O\left( \left\vert u\right\vert ^{\alpha +1}\right) $ \ for $%
u\rightarrow 0$, $\alpha \geq 0$ be a positive integer. If $u,$ $\upsilon
\in Y_{\infty }^{s,2},$ $\left\Vert u\right\Vert _{H^{s}}\leq M$, \ $%
\left\Vert \upsilon \right\Vert _{H^{s}}\leq M$ and $\left\Vert u\right\Vert
_{\infty }\leq M$, \ $\left\Vert \upsilon \right\Vert _{\infty }\leq M,$
then 
\[
\left\Vert f(u)-f(u)\right\Vert _{H^{s}}\leq C\left( M\right) \left[ \left(
\left\Vert u\right\Vert _{\infty }-\left\Vert \upsilon \right\Vert _{\infty
}\right) \left( \left\Vert u\right\Vert _{H^{s}}+\left\Vert \upsilon
\right\Vert _{H^{s}}\right) \right. 
\]%
\[
\left( \left\Vert u\right\Vert _{\infty }+\left\Vert \upsilon \right\Vert
_{\infty }\right) ^{\alpha -1}, 
\]%
\[
\left\Vert f(u)-f(\upsilon \right\Vert _{1}\leq C\left( M\right) \left(
\left\Vert u\right\Vert _{\infty }+\left\Vert \upsilon \right\Vert _{\infty
}\right) ^{\alpha -1}\left( \left\Vert u\right\Vert _{2}+\left\Vert \upsilon
\right\Vert _{2}\right) \left\Vert u-\upsilon \right\Vert _{2}. 
\]

By reasoning as in Theorem 3.1 and $\left[ \text{13, Theorem 1.1}\right] $
we have:

\ \textbf{Theorem 3.2. }Let the Condition 3.1 hold. Assume $f\in C^{k}\left(
R\right) ,$\ with $k$ an integer $k\geq s>\frac{n}{2},$ satisfies $f\left(
u\right) =O\left( \left\vert u\right\vert ^{\alpha +1}\right) $ \ for $%
u\rightarrow 0.$ Then there exists a constant $\delta >0,$ such that for any 
$\varphi ,$ $\psi \in Y_{1}^{s,2}$ satisfying 
\begin{equation}
\left\Vert \varphi \right\Vert _{Y_{1}^{s,2}}+\left\Vert \psi \right\Vert
_{Y_{1}^{s,2}}\leq \delta ,  \tag{3.19}
\end{equation}%
problem $\left( 1.1\right) -\left( 1.2\right) $ has a unique local strange
solution $u\in C^{\left( 2\right) }\left( \left[ 0,\right. \left. \infty
\right) ;Y^{s,2}\right) $. Moreover,

\begin{equation}
\sup_{0\leq t<\infty }\left( \left\Vert u\right\Vert _{Y^{s,2}}+\left\Vert
u_{t}\right\Vert _{Y^{s,2}}\right) \leq C\delta ,  \tag{3.20}
\end{equation}%
where the constant $C$ only depends on $f$ \ and initial data.

\textbf{Proof. \ }Consider a metric space \ defined by 
\[
W^{s}=\left\{ u\in ^{\left( 2\right) }\left( \left[ 0,\right. \left. \infty
\right) ;Y^{s,2}\right) ,\text{ }\left\Vert u\right\Vert _{W^{s}}\leq
3C_{0}\delta \right\} , 
\]%
equipped with the norm%
\[
\left\Vert u\right\Vert _{W^{s}}=\sup\limits_{t\geq 0}\left( \left\Vert
u\right\Vert _{Y_{\infty }^{s,2}}+\left\Vert u_{t}\right\Vert _{Y_{\infty
}^{s,2}}\right) , 
\]%
where $\delta >0$ \ satisfies $\left( 3.19\right) $ and $C_{0}$ \ is a
constant in Theorem 2.1. \ It is easy to prove that $W^{s}$ is a complete
metric space. From Sobolev imbedding theorem we know that $\left\Vert
u\right\Vert _{\infty }\leq 1$ if we take that $\delta $ is enough small.
Consider the problem $\left( 3.4\right) $. From Lemma 3.6 we get that $%
f\left( u\right) \in L^{2}\left( 0,T;Y_{1}^{s,2}\right) $ for any $T>0$.
Thus the problem $\left( 3.4\right) $ has a unique solution which can be
written as $\left( 3.5\right) .$ We should prove that the operator $G\left(
u\right) $\ defined by $\left( 3.5\right) $ is strictly contractive if $%
\delta $ is suitable small. \ In fact, by $(2.9)$ in Theorem 2.1 and Lemma
3.6 we get 
\[
\left\Vert G\left( u\right) \right\Vert _{\infty }+\left\Vert G_{t}\left(
u\right) \right\Vert _{\infty }\leq C_{0}\left[ \left\Vert \varphi
\right\Vert _{Y_{2}^{s,2}}\right. +\left\Vert \psi \right\Vert
_{Y_{2}^{s,2}}+ 
\]

\[
\left. \dint\limits_{0}^{t}\left\Vert K\left( u\right) \left( .,\tau \right)
\right\Vert _{Y_{2}^{s,2}}d\tau \right] \leq C_{0}\delta
+C_{0}\dint\limits_{0}^{t}\left\Vert K\left( u\right) \left( .,\tau \right)
\right\Vert _{Y_{2}^{s,2}}d\tau \leq 
\]%
\[
C_{0}\delta +C\dint\limits_{0}^{t}\left[ \left\Vert u\left( .,\tau \right)
\right\Vert _{2}^{2}\left\Vert u\left( .,\tau \right) \right\Vert _{\infty
}^{\alpha -1}+\left\Vert u\left( \tau \right) \right\Vert
_{H^{s}}^{2}\left\Vert u\left( \tau \right) \right\Vert _{\infty }^{\alpha }%
\right] d\tau \leq 
\]%
\begin{equation}
C_{0}\delta +C\left\Vert u\right\Vert _{W^{s}}^{\alpha +1}.  \tag{3.21}
\end{equation}%
On the other hand, by $(2.20)$ in Theorem 2.2 and Lemma 3.6 we have%
\[
\left\Vert G\left( u\right) \right\Vert _{H^{s}}+\left\Vert G_{t}\left(
u\right) \right\Vert _{H^{s}}\leq C_{0}\left[ \left\Vert \varphi \right\Vert
_{H^{s}}\right. +\left\Vert \psi \right\Vert _{H^{s}}+ 
\]

\[
\left. \dint\limits_{0}^{t}\left\Vert K\left( u\right) \left( .,\tau \right)
\right\Vert _{H^{s}}d\tau \right] \leq C_{0}\delta
+C_{0}\dint\limits_{0}^{t}\left\Vert K\left( u\right) \left( .,\tau \right)
\right\Vert _{H^{s}}d\tau \leq 
\]%
\begin{equation}
C_{0}\delta +C\dint\limits_{0}^{t}\left\Vert u\left( \tau \right)
\right\Vert _{H^{s}}\left\Vert u\left( \tau \right) \right\Vert _{\infty
}^{\alpha }d\tau \leq C_{0}\delta +C\left\Vert u\right\Vert _{W^{s}}^{\alpha
+1}.  \tag{3.22}
\end{equation}

\bigskip Therefore, combining $(3.21)$ with $(3.22)$ yields 
\begin{equation}
\left\Vert G\left( u\right) \right\Vert _{W^{s}}\leq 2C_{0}\delta
+C\left\Vert u\right\Vert _{W^{s}}^{\alpha +1}\text{.}  \tag{3.23}
\end{equation}

\bigskip Taking that $\delta $ is enough small such that $C\left(
3C_{0}\delta \right) ^{\alpha }$ $<1/3$, from $\left( 3.23\right) $ and from
Theorems 2.1, 2.2 we dedused that $G$ maps $W^{s}$ into $W^{s}$. Then, by
reasoning as the remaining part of $\left[ \text{13, Theorem 1.1}\right] $
we obtain that $G$ :$W^{s}\rightarrow W^{s}$ is strictly contractive. Using
the contraction mapping principle, we know that $G(u)$ has a unique fixed
point $u(x,t)\in $ $C^{2}([0,\infty ),H^{s})$ and $u(x,t)$ is the solution
of the problem $(1.1)-(1.2)$.

We claim that the solution $u(x,t)$ of the problem $(1.1)-(1.2)$ is also
unique in $C^{2}([0,\infty ),H^{s})$. In fact, let $u_{1}$ and $u_{2}$ be
two solutions of the problem $(1.1)-(1.2)$ and $u_{1}$, $u_{2}\in
C^{2}([0,\infty ),H^{s})$. Let $u=u_{1}-u_{2}$; then%
\[
u_{tt}-a\Delta u+b\ast u=\Delta \left[ g\ast \left( f\left( u_{1}\right)
-f\left( u_{2}\right) \right) \right] . 
\]

\bigskip This fact is derived in a similar way as in Theorem 3.2, by using
Theorems 2.1, 2.2 and Gronwall's inequality.

\textbf{Condition 3.2. }Let the Condition 2.1 holds. Assume $f\in C^{\left[ s%
\right] +1}\left( R\right) $ with $f(0)=0$ for some $s\geq 0$;

(2) Assume that the kernel $g$ is an integrable function whose Fourier
transform satisfies 
\[
0\leq \hat{g}\left( \xi \right) \leq \left( 1+\left\vert \xi \right\vert
^{2}\right) ^{-\frac{r}{2}}\text{ for all }\xi \in R^{n}\text{ and }r\geq 2. 
\]

\textbf{Theorem 3.3. }Let the Condition 3.2 hold. Moreover, $s\geq 0$ and $%
r\geq 2$. Then there is some $T>0$ such that the multipoint IVB $(1.1)-(1.2)$
is well posed with solution in $C^{2}\left( \left[ 0,T\right] ;H^{s}\right) $
for initial data $\varphi ,$ $\psi \in H^{s}.$

\textbf{Proof. }Consider the convolution operator $u\rightarrow Ku=\Delta %
\left[ g\ast f\left( u\right) \right] .$ In view of assumptios we have \ 
\begin{equation}
\left\Vert \Delta g\ast \upsilon \right\Vert _{H^{s}}\lesssim \left\Vert
\left( 1+\xi \right) ^{\frac{s}{2}}\left\vert \xi \right\vert ^{2}\hat{g}%
\left( \xi \right) \hat{\upsilon}\left( \xi \right) \right\Vert \lesssim
\left\Vert \upsilon \right\Vert _{H^{s}},  \tag{3.24}
\end{equation}%
i.e. $\Delta g\ast \upsilon $ is a bounded linear operator on $H^{s}.$ Then
by Corollary 3.1, $K\left( u\right) $ is locally Lipschitz on $H^{s}$. Then
by reasoning as in Theorem 3.2 and $\left[ \text{13, Theorem 1.1}\right] $
we obtain that $G$: $H^{s}\rightarrow H^{s}$ is strictly contractive. Using
the contraction mapping principle, we get that the operator $G(u)$ defined
by $\left( 3.5\right) $ has a unique fixed point $u(x,t)\in $ $%
C^{2}([0,\infty ),H^{s})$ \ and $u(x,t)$ is the solution of the problem $%
(1.1)-(1.2)$. Moreover, we show that the solution $u(x,t)$ of $(1.1)-(1.2)$
is also unique in $C^{2}([0,\infty ),H^{s})$. In fact, let $u_{1}$ and $%
u_{2} $ be two solutions of the problem $(1.1)-(1.2)$ and $u_{1}$, $u_{2}\in
C^{2}([0,\infty ),H^{s})$. Let $u=u_{1}-u_{2}$; then%
\[
u_{tt}-a\Delta u+b\ast u=\Delta \left[ g\ast \left( f\left( u_{1}\right)
-f\left( u_{2}\right) \right) \right] . 
\]

\bigskip This fact is derived in a similar way as in Theorem 3.2, by using
Theorems 2.1, 2.2 and Gronwall's inequality.

\textbf{Theorem 3.4. }Let the Condition 3.2 hold and $r>2+\frac{n}{2}$. Then
there is some $T>0$ such that the multipoint IVB $(1.1)-(1.2)$ is well posed
with solution in $C^{2}\left( \left[ 0,T\right] ;Y_{\infty }^{s,2}\right) $
for initial data $\varphi ,$ $\psi \in Y_{\infty }^{s,2}.$

\textbf{Proof. } All we need here, is to show that $K\ast f(u)$ is Lipschitz
on $Y_{\infty }^{s,2}$. Indeed, by reasoning as in Theorem 3.3 we have 
\[
\left\Vert \Delta g\ast \upsilon \right\Vert _{H^{s+r-2}}\lesssim \left\Vert
\left( 1+\xi \right) ^{\frac{s+r-2}{2}}\left\vert \xi \right\vert ^{2}\hat{g}%
\left( \xi \right) \hat{\upsilon}\left( \xi \right) \right\Vert \lesssim
\left\Vert \upsilon \right\Vert _{H^{s}}, 
\]

\bigskip Then $\Delta g\ast \upsilon $ is a bounded linear map from $H^{s}$
into $H^{s+r-2}$. Since $s\geq 0$ and $r$ $>2+\frac{n}{2}$ \ we get $s+r-2>%
\frac{n}{2}.$ The Sobolev embedding theorem implies that $\Delta g\ast
\upsilon $\ is a bounded linear map from $Y_{\infty }^{s,2}$ into $Y_{\infty
}^{s,2}$. Lemma 3.4 implies the Lipschitz condition on $Y_{\infty }^{s,2}$.
Then, by reasoning as in Theorem 3.3 we obtain the assertion.

The solution in theorems 3.2-3.4 can be extended to a maximal interval $%
[0,T_{\max }),$ where finite $T_{\max }$ is characterized by the blow-up
condition 
\[
\limsup\limits_{T\rightarrow T_{\max }}\left\Vert u\right\Vert _{Y_{\infty
}^{s,2}}=\infty . 
\]

\textbf{Lemma 3.8.} Suppose the conditions of theorems 3.4, 3.5 hold and $u$
is the solution of multipoint IVP $(1.1)-(1.2).$ Then there is a global
solution if for any $T<\infty $ we have%
\begin{equation}
\sup_{t\in \left[ 0,\right. \left. T\right) }\left( \left\Vert u\right\Vert
_{Y_{\infty }^{s,p}}+\left\Vert u_{t}\right\Vert _{Y_{\infty }^{s,p}}\right)
<\infty .  \tag{3.25}
\end{equation}

\textbf{Proof. }Indeed, by reasoning as in the second part of the proof of
Theorem 3.1, by using a continuation of local solution of $(1.1)-(1.2)$ and
assuming contrary that, $\left( 3.25\right) $ holds and $T_{0}<\infty $\ \
we obtain contradiction, i.e. we get $T_{0}=T_{\max }=\infty .$

\begin{center}
\textbf{4. Conservation of energy and global existence. }
\end{center}

In this section, we prove the existence and the uniqueness of the global
strong solution and the global classical solution for the problem $%
(1.1)-(1.2).$ \ For this purpose, we are going to make a priori estimates of
the local strong solution for the problem $(1.1)-(1.2).$

\textbf{Condition 4.1. }Let the Condition 2.1 holds.\textbf{\ }Assume that
the kernel $g$ is an integrable function whose Fourier transform satisfies 
\[
0<\hat{g}\left( \xi \right) \leq \left( 1+\left\vert \xi \right\vert
^{2}\right) ^{-\frac{r}{2}}\text{ for all }\xi \in R^{n}\text{ and }r\geq 2. 
\]

Let $F^{-1}$ denote the inverse Fourie rtransform. We consider the operator $%
B$ defined by 
\[
u\in D\left( B\right) =H^{s},\text{ }Bu=F^{-1}\left[ \left\vert \xi
\right\vert ^{-1}\left( \hat{g}\left( \xi \right) \right) ^{-\frac{-1}{2}}%
\hat{u}\left( \xi \right) \right] , 
\]

\bigskip Then it is clear to see that 
\begin{equation}
B^{-2}u=-\Delta g\ast u,\text{ }B^{-1}u=F^{-1}\left[ \left\vert \xi
\right\vert \left( \hat{g}\left( \xi \right) \right) ^{\frac{-1}{2}}\hat{u}%
\left( \xi \right) \right] .  \tag{4.1}
\end{equation}

First, we show the following\bigskip

\textbf{Lemma 4.1. }Suppose the conditions of theorems 3.4, 3.5 hold with $%
\hat{g}\left( \xi \right) >0$ and the solution of multipoint IVP $%
(1.1)-(1.2) $ exists in $C^{2}\left( \left[ 0,T\right] ;Y_{\infty
}^{s,2}\right) $ for some $s\geq 0.$ If $B\varphi \in L^{2}$ and $B\psi \in
L^{2},$ then $Bu,$ $Bu_{t}\in C^{1}\left( \left[ 0,\right. \left. T\right)
;L^{2}\right) .$

\bigskip \textbf{Proof.} By Lemma 2.1, problem $\left( 1.1\right) -\left(
1.2\right) $ is equ\i valent to following integral equation ,%
\begin{equation}
u\left( x,t\right) =\left[ S_{1}\left( x,t\right) \varphi +S_{2}\left(
x,t\right) \psi +\Phi \left( g\ast f\left( u\right) \right) \right] + 
\tag{4.2}
\end{equation}%
\[
\frac{1}{2\eta }\dint\limits_{0}^{t}F^{-1}\left[ \sin \eta \left( t-\tau
\right) \left\vert \xi \right\vert ^{2}\hat{g}\left( \xi \right) \hat{f}%
\left( G\left( u\right) \left( \xi \right) \right) \right] d\tau , 
\]%
where $S_{1}\left( x,t\right) $, $S_{2}\left( x,t\right) ,$ $\Phi $ are are
operator functions defined by $\left( 2.10\right) $ and $\left( 2.11\right) $%
, $\ $where $g$ replaced by $g\ast f\left( u\right) .$

From $\left( 4.2\right) $ we get that 
\[
u_{t}\left( x,t\right) =\left[ \frac{d}{dt}S_{1}\left( x,t\right) \varphi +%
\frac{d}{dt}S_{2}\left( x,t\right) \psi +\frac{d}{dt}\Phi \left( g\ast
f\left( u\right) \right) \right] + 
\]%
\begin{equation}
\frac{1}{2}\dint\limits_{0}^{t}F^{-1}Q\left( \xi ,t-\tau \right) d\tau , 
\tag{4.3}
\end{equation}%
where 
\[
Q\left( \xi ,t-\tau \right) =\cos \eta \left( t-\tau \right) \left\vert \xi
\right\vert ^{2}\hat{g}\left( \xi \right) \hat{f}\left( u\right) \left( \xi
\right) , 
\]%
\[
S_{1}\left( x,t\right) \varphi =F^{-1}\left\{ D_{0}^{-1}\left( \xi \right) 
\left[ \left( 1-\dsum\limits_{k=1}^{m}\beta _{k}\cos \varkappa _{k}\right)
\sin \left( \eta t\right) \right] \right. ,+\text{ } 
\]%
\begin{equation}
\left[ \eta ^{-1}\left( \xi \right) \dsum\limits_{k=1}^{m}\beta _{k}\sin
\varkappa _{k}\cos \left( \eta t\right) \right] \left. \hat{\varphi}\left(
\xi \right) \right\} ,  \tag{4.4}
\end{equation}

\[
\frac{d}{dt}S_{1}\left( x,t\right) \varphi =F^{-1}\left\{ D_{0}^{-1}\left(
\xi \right) \eta \left[ \left( 1-\dsum\limits_{k=1}^{m}\beta _{k}\cos
\varkappa _{k}\right) \cos \left( \eta t\right) \right] \right. -\text{ } 
\]%
\[
\left[ \dsum\limits_{k=1}^{m}\beta _{k}\sin \varkappa _{k}\sin \left( \eta
t\right) \right] \left. \hat{\varphi}\left( \xi \right) \right\} , 
\]%
\[
S_{2}\left( x,t\right) \psi =F^{-1}\left\{ \left[ \eta ^{-1}\left( \xi
\right) D_{0}^{-1}\left( \xi \right) \left( \dsum\limits_{k=1}^{m}\alpha
_{k}\sin \varkappa _{k}\right) \sin \left( \eta t\right) \right] \right. +%
\text{ } 
\]

\begin{equation}
\eta ^{-1}\left( \xi \right) D_{0}^{-1}\left( \xi \right) \left. \left[
\left( 1-\dsum\limits_{k=1}^{m}\alpha _{k}\cos \varkappa _{k}\right) \cos
\left( \eta t\right) \right] \right\} \hat{\psi}\left( \xi \right) , 
\tag{4.5}
\end{equation}%
\[
\frac{d}{dt}S_{2}\left( x,t\right) \psi =F^{-1}\left\{ \left[
D_{0}^{-1}\left( \xi \right) \left( \dsum\limits_{k=1}^{m}\alpha _{k}\sin
\varkappa _{k}\right) \cos \left( \eta t\right) \right] \right. -\text{ } 
\]

\[
D_{0}^{-1}\left( \xi \right) \left. \left[ \left(
1-\dsum\limits_{k=1}^{m}\alpha _{k}\cos \varkappa _{k}\right) \sin \left(
\eta t\right) \right] \right\} \hat{\psi}\left( \xi \right) , 
\]

\begin{equation}
\Phi \left( x;t\right) =\Phi \left( g\ast f\left( u\right) \right)
=\dsum\limits_{j=}^{4}\dsum\limits_{k=1}^{m}F^{-1}\Phi _{jk}\left(
\left\vert \xi \right\vert ^{2}\hat{g}\left( \xi \right) \hat{f}\left(
u\right) \right) \left( \xi ;t\right) ,\text{ }  \tag{4.6}
\end{equation}

\bigskip Since 
\[
D_{0}^{-1}\left( \xi \right) ,\text{ }\left( 1-\dsum\limits_{k=1}^{m}\beta
_{k}\cos \varkappa _{k}\right) \sin \left( \eta t\right) ,\text{ }\eta
^{-1}\left( \xi \right) \dsum\limits_{k=1}^{m}\beta _{k}\sin \varkappa
_{k}\cos \left( \eta t\right) 
\]%
are uniformly bounded for fixet $t$ \ by $\left( 4.1\right) $, $\left(
4.4\right) $ $\left( 4.5\right) $ we have%
\begin{equation}
\left\Vert BS_{1}\left( x,t\right) \varphi \right\Vert _{L^{2}}=\left\Vert
F^{-1}\left[ \left\vert \xi \right\vert ^{-1}\left( \hat{g}\left( \xi
\right) \right) ^{-\frac{-1}{2}}\hat{u}\left( \xi \right) S_{1}\left(
x,t\right) \varphi \right] \right\Vert _{L^{2}}\lesssim \left\Vert \varphi
\right\Vert _{H^{s}}<\infty ,  \tag{4.7}
\end{equation}%
\[
\left\Vert BS_{2}\left( x,t\right) \varphi \right\Vert _{L^{2}}=\left\Vert
F^{-1}\left[ \left\vert \xi \right\vert ^{-1}\left( \hat{g}\left( \xi
\right) \right) ^{-\frac{-1}{2}}\hat{u}\left( \xi \right) S_{2}\left(
x,t\right) \psi \right] \right\Vert _{L^{2}}\lesssim \left\Vert \psi
\right\Vert _{H^{s}}<\infty . 
\]

For fixed $t$, we have $f(u)\in H^{s}.$ Since \ $D_{0}^{-1}\left( \xi
\right) ,$ $\cos \left( \eta t\right) ,$ $\sin \left( \eta t\right) $ are
uniformly bounded,\ from $\left( 4.1\right) ,$ $\left( 2.11\right) $ and $%
\left( 4.6\right) $\ we get 
\begin{equation}
\left\Vert B\Phi \right\Vert _{L^{2}}\leq \left\Vert F^{-1}\left[ \left\vert
\xi \right\vert ^{-1}\left( \hat{g}\left( \xi \right) \right) ^{-\frac{-1}{2}%
}\left\vert \xi \right\vert ^{2}\hat{g}\left( \xi \right) \hat{f}\left(
u\right) \right] \right\Vert _{L^{2}}\lesssim \left\Vert f\left( u\right)
\right\Vert _{H^{s}}<\infty .  \tag{4.8}
\end{equation}

\bigskip From $\left( 4.8\right) $ we have 
\begin{equation}
\left\Vert B\Phi \right\Vert _{L^{2}}\leq \left\Vert F^{-1}\left[ \left\vert
\xi \right\vert ^{-1}\left( \hat{g}\left( \xi \right) \right) ^{-\frac{-1}{2}%
}Q\left( \xi ,t-\tau \right) \right] \right\Vert _{L^{2}}\leq  \tag{4.9}
\end{equation}%
\[
\left\Vert F^{-1}\left[ \left\vert \xi \right\vert ^{-1}\left( \hat{g}\left(
\xi \right) \right) ^{-\frac{-1}{2}}\left\vert \xi \right\vert ^{2}\hat{g}%
\left( \xi \right) \hat{f}\left( u\right) \right] \right\Vert
_{L^{2}}\lesssim \left\Vert f\left( u\right) \right\Vert _{H^{s}}<\infty . 
\]%
Then from $\left( 4.2\right) ,$ $\left( 4.7\right) ,$ $\left( 4.8\right) $
and $\left( 4.9\right) $ we obtain the assertion.

\textbf{Remark 4.1.} Due to nonlocality of initial conditions the additional
conditions appears in Theorem 4.1. For classical Cauchy problem this extra
conditions are not required

\textbf{Lemma 4.2.} Assume the conditions of theorems 3.4, 3.5 hold with $%
a=0,$ $\hat{g}\left( \xi \right) >0$ and%
\[
\hat{b}\left( \xi \right) =O\left( 1+\left\vert \xi \right\vert ^{2}\right)
^{\frac{s-r}{2}+1}.
\]%
Suppose the solution of $(1.1)-(1.2)$ exists in $C^{2}\left( \left[
0,\right. \left. T\right) ;Y_{\infty }^{s,2}\right) $ for some $s\geq 0.$ If 
$B\psi \in L^{2},$ $Bu_{t}\left( x,\lambda _{k}\right) \in L^{2},k=1,2,...,m,
$ then $Bu_{t}\in C^{2}\left( \left[ 0,\right. \left. T\right) ;L^{2}\right)
.$ Moreover, if $B\varphi \in L^{2},$ $Bu\left( x,\lambda _{k}\right) \in
L^{2}$ then $Bu\in C^{1}\left( \left[ 0,\right. \left. T\right)
;L^{2}\right) .$

\textbf{Proof. }\ Integrating the equation $\left( 1.1\right) $ for $a=0,$
twice and calculating the resulting double integral as an iterated integral,
we have 
\[
u\left( x,t\right) =\varphi \left( x\right) +\dsum\limits_{k=1}^{m}\alpha
_{k}u\left( x,\lambda _{k}\right) +t\left[ \psi \left( x\right)
+\dsum\limits_{k=1}^{m}\beta _{k}u_{t}\left( x,\lambda _{k}\right) \right] - 
\]%
\begin{equation}
\dint\limits_{0}^{t}\left( t-\tau \right) \left( b\ast u\right) \left(
x,\tau \right) d\tau +\dint\limits_{0}^{t}\left( t-\tau \right) \Delta
\left( g\ast f\left( u\right) \right) \left( x,\tau \right) d\tau , 
\tag{4.10}
\end{equation}

\[
u_{t}\left( x,t\right) =\psi \left( x\right) +\dsum\limits_{k=1}^{m}\beta
_{k}u_{t}\left( x,\lambda _{k}\right) - 
\]%
\begin{equation}
\dint\limits_{0}^{t}\left( b\ast u\right) \left( x,\tau \right) d\tau
+\dint\limits_{0}^{t}\Delta \left( g\ast f\left( u\right) \right) \left(
x,\tau \right) d\tau .  \tag{4.11}
\end{equation}

\bigskip From $\left( 4.1\right) $ and $\left( 4.11\right) $ for fixed $t$
and $\tau $ we get 
\[
\left\Vert Bu_{t}\left( x,t\right) \right\Vert _{L^{2}}=\left\Vert B\psi
\left( x\right) \right\Vert _{L^{2}}+\dsum\limits_{k=1}^{m}\beta
_{k}\left\Vert Bu_{t}\left( x,\lambda _{k}\right) \right\Vert _{L^{2}}- 
\]%
\begin{equation}
\dint\limits_{0}^{t}\left\Vert B\left( b\ast u\right) \left( x,\tau \right)
\right\Vert _{L^{2}}d\tau -\dint\limits_{0}^{t}\left\Vert B^{-1}f\left(
u\right) \left( x,\tau \right) \right\Vert _{L^{2}}d\tau .  \tag{4.12}
\end{equation}

\bigskip

\bigskip By assumption on $b$, $g$ and by $\left( 4.1\right) $ for fixed $%
\tau $ we have%
\begin{equation}
\left\Vert B\left( b\ast u\right) \left( x,\tau \right) \right\Vert
_{L^{2}}\leq \left\Vert F^{-1}\left[ \left\vert \xi \right\vert ^{-1}\hat{b}%
\left( \xi \right) \left( \hat{g}\left( \xi \right) \right) ^{-\frac{-1}{2}}%
\hat{u}\left( \xi ,\tau \right) \right] \right\Vert _{L^{2}}\lesssim
\left\Vert u\left( .,\tau \right) \right\Vert _{H^{s}}.  \tag{4.13}
\end{equation}

\bigskip Moreover, by Lemma 3.3 for all $t$ we have $f\left( u\right) \in
H^{s}.$ Also 
\begin{equation}
\left\Vert B\left( f\left( u\right) \right) \left( x,\tau \right)
\right\Vert _{L^{2}}\lesssim \left\Vert f\left( u\right) \left( .,\tau
\right) \right\Vert _{H^{s}}.  \tag{4.14}
\end{equation}

\bigskip Then from $\left( 4.12\right) -\left( 4.14\right) $ we obtain $%
Bu_{t}\in C^{2}\left( \left[ 0,\right. \left. T\right) ;L^{2}\right) .$ The
second statement follows similarly from $\left( 4.10\right) .$

From Lemma 4.2 we obtain the following result.

\textbf{Result 4.1.} Assume the conditions of theorems 3.4, 3.5 hold with $%
a=0,$ $\hat{g}\left( \xi \right) >0,$ $\alpha _{k}=\beta _{k}=0$ and%
\[
\hat{b}\left( \xi \right) =O\left( 1+\left\vert \xi \right\vert ^{2}\right)
^{\frac{s-r}{2}+1}.
\]%
Suppose the solution of $(1.1)-(1.2)$ exists in $C^{2}\left( \left[ 0,T%
\right] ;Y_{\infty }^{s,2}\right) $ for some $s\geq 0.$ If $B\psi \in L^{2}$
then $Bu_{t}\in C^{2}\left( \left[ 0,\right. \left. T\right) ;L^{2}\right) .$
Moreover, if $B\varphi \in L^{2},$ then $Bu\in C^{1}\left( \left[ 0,\right.
\left. T\right) ;L^{2}\right) .$

Here, 
\[
G\left( \tau \right) =\dint\limits_{0}^{\tau }g\left( s\right) ds.
\]

\bigskip \textbf{Lemma 4.3. }Assume the all conditions of Lemma 4.2 are
satisfied. Let $B\psi \in L^{2},$ $Bu_{t}\left( x,\lambda _{k}\right) \in
L^{2},$ $k=1,2,...,m$ and $G\left( \varphi \right) \in L^{1}$. Then for any $%
t\in \left[ 0,\right. \left. T\right) $ the energy 
\begin{equation}
E\left( t\right) =\left\Vert Bu_{t}\right\Vert _{L^{2}}^{2}+a\left\Vert
F^{-1}\hat{g}\ast u\right\Vert _{L^{2}}^{2}+\left\Vert B\left( b\ast
u\right) \right\Vert _{L^{2}}^{2}+2\dint\limits_{R^{n}}G\left( u\right) dx 
\tag{4.15}
\end{equation}%
is constant $\left[ 0,\right. \left. T\right) .$ 

\textbf{Proof. } By use of \ equation $\left( 1.1\right) $, it follows from
straightforward calculation that 
\[
\frac{d}{dt}E\left( t\right) =2\left( Bu_{tt},Bu_{t}\right) +2a\left( F^{-1}%
\hat{g}\ast u,\left( F^{-1}\hat{g}\ast u\right) _{t}\right) +
\]

\[
2\left[ B\left( b\ast u\right) ,B\left( b\ast u\right) u_{t}\left( t\right) %
\right] +2\left( f\left( u\right) ,u_{t}\right) =2\left(
B^{2}u_{tt},u_{t}\right) +
\]%
\[
2\left( B^{2}\left( b\ast u\right) ,\left( b\ast u\right) u_{t}\left(
t\right) \right) +2\left[ B^{2}\left( b\ast u\right) ,\left( b\ast u\right)
u_{t}\left( t\right) \right] =
\]%
\[
2B^{2}\left[ \left( u_{tt}-a\Delta u+b\ast u+\Delta \left[ g\ast f\left(
u\right) \right] ,u_{t}\right) \right] =0,
\]%
where $\left( u,\upsilon \right) $  denotes the inner product of $L^{2}$
space. \ Integrating the above equality with respect to $t$, we have $\left(
4.15\right) $.

\bigskip By using the above lemmas we obtain the following results

\textbf{Theorem 4.1. }Let the Condition 3.2 hold $a=0,$ $\hat{g}\left( \xi
\right) >0$ and%
\[
\hat{b}\left( \xi \right) =O\left( 1+\left\vert \xi \right\vert ^{2}\right)
^{\frac{s-r}{2}+1}.
\]
Moreover, let $B\psi \in L^{2},$ $Bu_{t}\left( x,\lambda _{k}\right) \in
L^{2},$ $k=1,2,...,m$ and $G\left( \varphi \right) \in L^{1}$, $s\geq 0$ , $%
r>3$ and there is some $k>0$ so that \ $G\left( r\right) \geq -kr^{2}$ for $%
r\in \mathbb{R}$. Then there is some $T>0$ such that the multipoint IVB $%
(1.1)-(1.2)$ has a global solution $u\in C^{2}\left( \left[ 0,\right. \left.
\infty \right) ;Y_{\infty }^{s,2}\right) .$

\bigskip

\textbf{References}

\begin{quote}
\ \ \ \ \ \ \ \ \ \ \ \ \ \ \ \ \ \ \ \ \ \ \ \ 
\end{quote}

\begin{enumerate}
\item A. C. Eringen, Nonlocal Continuum Field Theories (2002), New York,
Springer.

\item Z. Huang, Formulations of nonlocal continuum mechanics based on a new
definition of stress tensor Acta Mech. (2006)187, 11--27.

\item C. Polizzotto, Nonlocal elasticity and related variational principles
Int. J. Solids Struct. ( 2001) 38 7359--80.

\item C. A. Silling, Reformulation of elasticity theory for discontinuities
and long-range forces J. Mech. Phys. Solids (2000)48 175-209.

\item M. Arndt and M. Griebel, Derivation of higher order gradient continuum
models from atomistic models for crystalline solids Multiscale Modeling
Simul. (2005)4, 531--62..

\item X. Blanc, C. LeBris, P. L. Lions, Atomistic to continuum limits for
computational materials science, ESAIM--- Math. Modelling Numer. Anal.
(2007)41, 391--426.

\item A. De Godefroy, Blow up of solutions of a generalized Boussinesq
equation IMA J. Appl. Math.(1998) 60 123--38.

\item A. Constantin and L. Molinet, The initial value problem for a
generalized Boussinesq equation, Diff.Integral Eqns. (2002)15, 1061--72.

\item G. Chen and S. Wang, Existence and nonexistence of global solutions
for the generalized IMBq equation Nonlinear Anal.---Theory Methods Appl.
(1999)36, 961--80.

\item M. Lazar, G. A. Maugin and E. C. Aifantis, On a theory of nonlocal
elasticity of bi-Helmholtz type and some applications Int. J. Solids and
Struct. (2006)43, 1404--21.

\item N. Duruk, H.A. Erbay and A. Erkip, Global existence and blow-up for a
class of nonlocal nonlinear Cauchy problems arising in elasticity,
Nonlinearity, (2010)23, 107--118.

\item S. Wang, G. Chen, Small amplitude solutions of the generalized IMBq
equation, J. Math. Anal. Appl. 274 (2002) 846--866.

\item Wang S and Chen G, Cauchy problem of the generalized double dispersion
equation Nonlinear Anal.--- Theory Methods Appl. (2006 )64 159--73.

\item J.L. Bona, R.L. Sachs, Global existence of smooth solutions and
stability of solitary waves for a generalized Boussinesq equation, Comm.
Math. Phys. 118 (1988), 15--29.

\item F. Linares, Global existence of small solutions for a generalized
Boussinesq equation, J. Differential Equations 106 (1993), 257--293.

\item Y. Liu, Instability and blow-up of solutions to a generalized
Boussinesq equation, SIAM J. Math. Anal. 26 (1995), 1527--1546.

\item V.G. Makhankov, Dynamics of classical solutions (in non-integrable
systems), Phys. Lett. C 35(1978), 1--128.

\item G.B. Whitham, Linear and Nonlinear Waves, Wiley--Interscience, New
York, 1975.

\item N.J. Zabusky, Nonlinear Partial Differential Equations, Academic
Press, New York, 1967.

\item A. Ashyralyev, N. Aggez, Nonlocal boundary value hyperbolic problems
involving Integral conditions, Bound.Value Probl., 2014 V (2014):214.

\item L. S. Pulkina, A non local problem with integral conditions for
hyperbolice quations, Electron.J.Differ.Equ.(1999)45, 1-6.

\item Girardi, M., Lutz, W., Operator-valued Fourier multiplier theorems on $%
L_{p}(X)$ and geometry of Banach spaces, J. Funct. Anal., ( 2003)204(2),
320--354.

\item H. Triebel, Interpolation theory, Function spaces, Differential
operators, North-Holland, Amsterdam, 1978.

\item H. Triebel, Fractals and spectra, Birkhauser Verlag, Related to
Fourier analysis and function spaces, Basel, 1997.

\item L. Nirenberg, On elliptic partial differential equations, Ann. Scuola
Norm. Sup. Pisa (1959)13 , 115--162.

\item S. Klainerman, Global existence for nonlinear wave equations, Comm.
Pure Appl. Math.(1980)33 , 43--101.

\item R. Coifman and Y. Meyer, Wavelets. Calder%
\'{}%
on-Zygmund and Multilinear Operators, Cambridge University Press, 1997.

\item T. Kato, G. Ponce, Commutator estimates and the Euler and
Navier--Stokes equations, Comm. Pure Appl. Math. (1988)41, 891--907.
\end{enumerate}

\end{document}